\documentclass[AMA,STIX1COL,doublespace]{WileyNJD-v2}
\articletype{Research Article}

\usepackage{amsmath,amsfonts,amsthm}
\usepackage{epsfig}
\usepackage{subfigure}
\usepackage{graphicx}
\usepackage{fancybox}
\usepackage{bm}
\usepackage{verbatim}
\usepackage{cancel}
\usepackage{color,colordvi}
\usepackage{enumerate}
\renewcommand{\Re}{{\mathbb{R}}}

\newcommand{\bsym}[1]{\boldsymbol{#1}}

\graphicspath{ {./figs/} }

\begin{document}
	
\title{Axisymmetric Virtual Elements For Problems of Elasticity and Plasticity}
	
\author[1]{L. L. Yaw$\mbox{}^{1*,}$}
%\author[2]{John J. John}
	
\authormark{L. L. Yaw}% \sc{and} John J. John}
	
\address[1]{\orgdiv{Engineering Department},
	\orgname{Walla Walla University},
	\orgaddress{College Place, WA 99324, USA}}
%\address[2]{\orgname{J Lab},
%	\orgaddress{Albuquerque, NM 87185}, \country{USA}}
\corres{$\mbox{}^*$Louie L. Yaw, Engineering Department,
	Walla Walla University, 100 SW 4th St, College Place, 
	WA 99324, USA \\
	\email{louie.yaw@wallawalla.edu}}

\abstract{
The virtual element method (VEM) allows discretization of elasticity and plasticity problems with polygons in 2D and polyhedrals in 3D.  The polygons (and polyhedrals) can have an arbitrary number of sides and can be concave or convex.  These features, among others, are attractive for meshing complex geometries.  However, to the author's knowledge axisymmetric virtual elements have not appeared before in the literature.  Hence, in this work a novel first order consistent axisymmetric virtual element method is applied to problems of elasticity and plasticity.  The VEM specific implementation details and adjustments needed to solve axisymmetric simulations are presented.  Representative benchmark problems including pressure vessels and circular plates are illustrated.  Examples also show that problems of near incompressibility are solved successfully.  Consequently, this research demonstrates that the axisymmetric VEM formulation successfully solves certain classes of solid mechanics problems.  The work concludes with a discussion of results for the current formulation and future research directions.}

\keywords{virtual elements; axisymmetric; nonlinear; elasticity; plasticity; near incompressibility; pressure vessels; circular plates}

\maketitle

% Section 1
\section{Introduction}\label{sec:intro}
Since its inception \cite{Beirao,Beirao2} the virtual element method (VEM) has attracted much attention.  Its use in solid mechanics for problems of elasticity \cite{Beirao:2013:VEL,Artioli} and plasticity \cite{TaylorVEMplast,Artioli2,Beirao3} has been explored by a variety of researchers for 2D and 3D problems  Yet, an axisymmetric VEM formulation is not yet present in the literature.  For certain classes of problems it is more computationally efficient to have an axisymmetric method available.  In this work a novel first order consistent axisymmetric virtual element method is applied to axisymmetric elasto-static problems.  Furthermore, small strain plasticity problems are solved with the proposed formulation as well.
The structure of the remainder of this paper follows.  In section \ref{sec:reviewVEM}, a brief review of virtual elements for 2D elasticity is provided with the reader directed to relevant references for more details.  VEM notation, used herein, is also provided.  The VEM review is followed by section \ref{sec:MVC}, which provides details on mean value coordinates (MVC).  Mean value coordinates are needed to modify the 2D VEM formulation into an axisymmetric formulation.  With the preceding sections in hand, section \ref{sec:axi} sets forth the necessary details to construct an axisymmetric formulation using VEM.  Section \ref{sec:elast} discusses details associated with elasto-statics problems.  In section \ref{sec:plast}, including plasticity in the formulation is presented.  Numerical implementation details, particular to this work, are discussed in section \ref{sec:numimplementation}.  Representative results of numerical simulations are provided in section \ref{sec:numresults}.  Simulation results are compared to benchmark problems or to known theoretical solutions.  Last, main findings and conclusions are provided with thoughts on improvements and future research directions in section \ref{sec:conclusions}.

%2 Review of 2D Elasticity VEM \ref{sec:reviewVEM}
%3 Mean value coordinates\ref{sec:MVC}
%4 Review of axisymmetric formulation\ref{sec:axi}
%5 Axisymm VEM elasticity \ref{sec:elast}
%6 Axisymm VEM plasticity \ref{sec:plast}
%7 Numerical implementation \ref{sec:numimplementation}
%8 Numerical results \ref{sec:numresults}
%9 Conclusions \ref{sec:conclusions}

% Section 2
\section{Review of Virtual Element Method (VEM) For 2D Elasticity}\label{sec:reviewVEM}
In 2013 \cite{Beirao} the virtual element method (VEM) appeared in the literature.  It is an attractive method for solving partial differential equations in science, mathematics, and engineering.  Due to similarities with the finite element method (FEM), VEM has been applied to a variety of engineering problems.  However, VEM domain discretizations allow for arbitrary polygons, rather than restriction to triangles and quadrilaterals.  The polygons need not all be the same and are permitted to be concave or convex.  These traits are appealing since they simplify mesh construction.  Polynomial consistency, that is linear, quadratic, or higher, is accomplished in the VEM formulation.  Furthermore, VEM easily allows for non-conforming discretizations (see Mengolini \cite{mengolini}).  The current work is confined to linear order ($k=1$) polynomial interpolation (consistency), and a 2D linear elasticity formulation is described in this section.    The goal here is to privde only minimum implementation details.  For additional information, interested researchers should consult the provided references.  With some exceptions, the presentation follows notational conventions provided in references Mengolini\cite{mengolini} and Yaw\cite{yawVEM}.

VEM is described by Sukumar and Tupek \cite{suku:elastodyn}.  In VEM, ``the basis functions are defined as the solution of a local elliptic partial differential equation", yet the basis functions are not calculated to utilize the method.  The VEM basis function are defined in a convenient way, but their actual form is not known.  As a result they are described as being \emph{virtual}, and the finite
element space for VEM as a \emph{virtual element space}.  Similar to FEM, by piecing together a local discretization space  $\boldsymbol{\mathcal{V}}_k(E)$ (for element E of order $k$) a global conforming space $\boldsymbol{\mathcal{V}}^h$ is constructed.  Yet, $\boldsymbol{\mathcal{V}}_k(E)$ contains polynomial trial and test functions that are $k$th order or less and may also contain nonpolynomial functions.  These traits are different than FEM.  For the weak form, the stiffness (bilinear form) and forces (linear form) are formed from elliptic polynomial projections of the VEM basis functions.  The computed projections are accomplished from the degrees of freedom for the linear elasticity problem and provide no additional approximation error.  The result is a stiffness that has two parts:  a consistency part tied to the chosen polynomial space, and a stability part that provides coercivity or invertibility.  Then, like FEM, an element by element assembly process is utilized to construct the global system of matrices.  VEM is analogous to a stabilized hourglass control finite element method~\cite{Flanagan:1981:AUS, Russo} which uses convex and or nonconvex polytopes.\cite{Cangiani:2015:HSV}

\subsection{The Continuous 2D Linear Elasticity Problem}
VEM is employed to solve 2D elasticity problems.
Recall, for elasticity problems, the standard weak form:  Find $\mathbf{u} \in \boldsymbol{\mathcal{V}}$ such that
\begin{equation}\label{E1}
	a(\mathbf{u},\mathbf{v})=L(\mathbf{v}) \quad \forall \mathbf{v}\in \boldsymbol{\mathcal{V}},
\end{equation}
where the bilinear form is
\begin{equation}\label{E2}
	a(\mathbf{u},\mathbf{v})=\int _{\Omega} \! \boldsymbol{\sigma}(\mathbf{u}):\boldsymbol{\varepsilon}(\mathbf{v}) \, \mathrm{d}\Omega,
\end{equation}
and the linear form is
\begin{equation}\label{E3}
	L(\mathbf{v})=\int_{\Omega} \! \mathbf{v}\cdot\mathbf{f} \, \mathrm{d}\Omega+\int_{\partial\Omega_t} \! \mathbf{v}\cdot\bar{\mathbf{t}} ,\ \mathrm{d}\partial\Omega.
\end{equation}\\

\noindent
\textbf{Remarks}
\begin{enumerate}[(i)]
	\item The space of functions $\boldsymbol{\mathcal{V}}$ is vector-valued containing components $v_x$ and $v_y$.  These functions are found in  first-order Sobolev space $\mathcal{H}^1(\Omega)$ and on the displacement boundary their value is zero.
	\item Trial functions are , $\mathbf{u} \in \boldsymbol{\mathcal{V}}$, and weight functions are, $\mathbf{v} \in \boldsymbol{\mathcal{V}}$.
	\item Throughout this work, unless otherwise noted, Voigt notation is intended in equations with column vectors and matrices.  One exception is \eqref{E2}, where tensor notation is used.
\end{enumerate}

\subsection{Discretization of the problem domain}
\begin{figure}
\centering
\subfigure[]{\epsfig{file =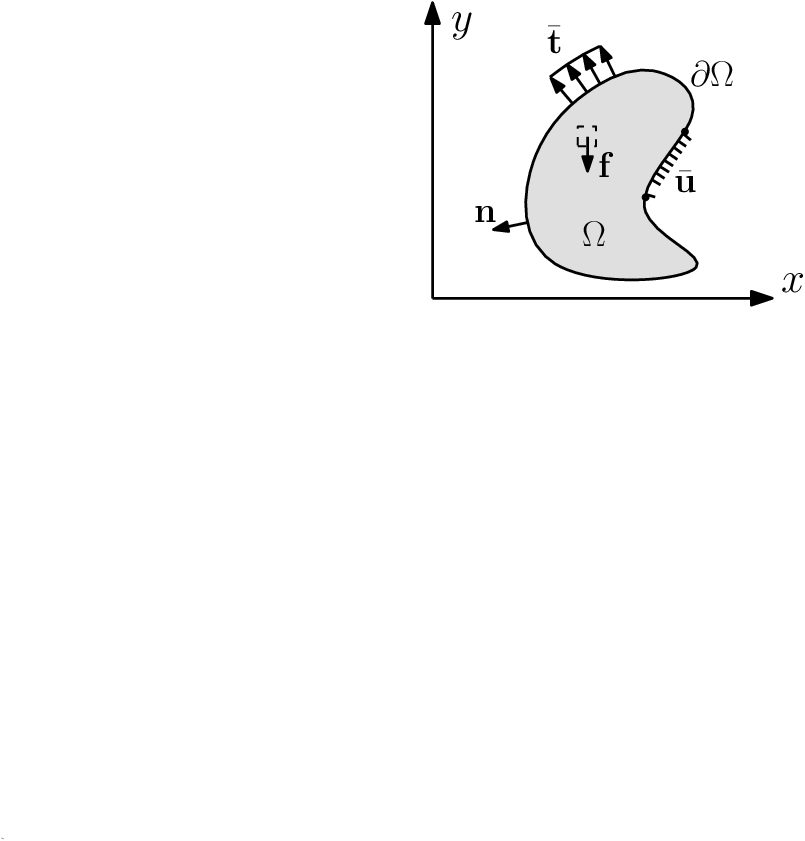,width=0.38\textwidth}}
\hspace*{0.4in}
\subfigure[]{\epsfig{file = 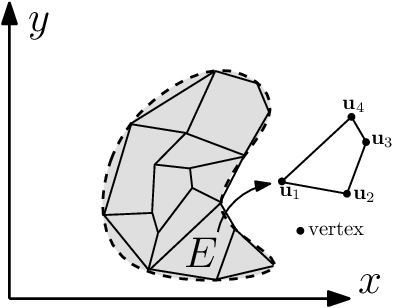,width=0.4\textwidth}}
\caption{2D Solid Domain: (a)~Elasticity problem with boundary conditions, (b)~Virtual element method domain discretization and example polygonal element with vector of nodal displacements labeling each vertex.}
\label{fig1}
\end{figure}

An example geometric domain of interest is provided in Figure~\ref{fig1}a.  Calculating field variables such as stress, strain, and displacements are the objective.  By discretizing the domain with polygons, Figure~\ref{fig1}b, it is then possible to use VEM.  For VEM the polygons can be convex or non-convex with arbitrary number of sides.  The interpolation space is necessarily populated with polynomial as well as non-polynomial functions as a consequence of selecting polygonal elements for discretization.  The resulting polygon elements are able to compatibly interface with adjacent elements along their sides.  Essentially, the polygon edges have polynomial interpolation functions, but polynomial plus (possibly) non-polynomial functions exist on the element interior.  Nevertheless, it is only necessary to know polynomial functions along the edges to implement VEM.  Consequently, the edge polynomials are chosen by selecting a polynomial order, $k$.  The present work chooses first order polynomials exclusively.
		
%Typical double figure

On element E, scalar-valued polynomials of order equal to $k$ or less are in the space $\mathcal{P}_k(E)$.  Consequently, $\boldsymbol{\mathcal{P}}_k\equiv [\mathcal{P}_k]^2$ is used to represent a 2D vector space of polynomials with two variables.  The expression $\mathbf{P}_k=\{\mathbf{p}_{\alpha}\}_{\alpha=1,...,n_k}$ denotes a basis for the polynomial space.  Illustrated below is a polynomial basis of order $k=1$.
\begin{equation}\label{pf1}
	\mathbf{P}_1=[\mathbf{p}_1,\ \mathbf{p}_2,\ \mathbf{p}_3,\ \mathbf{p}_4,\ \mathbf{p}_5,\ \mathbf{p}_6]
\end{equation}
or
\begin{equation}\label{pf2}
	\mathbf{P}_1=\left[\left ( \begin{array}{c}
	1\\
	0\\
	\end{array}\right ),
	\ \left ( \begin{array}{c}
	0\\
	1\\
	\end{array}\right ),
	\ \left ( \begin{array}{c}
	-\eta\\
	\xi\\
	\end{array}\right ),
	\ \left ( \begin{array}{c}
	\eta\\
	\xi\\
	\end{array}\right ),
	\ \left ( \begin{array}{c}
	\xi\\
	0\\
	\end{array}\right ),
	\ \left ( \begin{array}{c}
	0\\
	\eta\\
	\end{array}\right )
	\right].
\end{equation}
In equation \eqref{pf2}, scaled monomials are used to create the components.  The scaled monomials are defined as
\begin{equation}
	\xi=\left(\frac{x-\bar{x}}{h_E}\right), \quad \eta=\left(\frac{y-\bar{y}}{h_E}\right),
\end{equation}
where $\bar{\mathbf{x}}=(\bar{x},\bar{y})$ is the element $E$ centroid location and the diamter is $h_E$(i.e., diameter of smallest circle that encloses all element vertices).\\

\noindent
\textbf{Remarks}
\begin{enumerate}[(i)]
	\item It is possible to choose polynomials of higher order (see \cite{mengolini}).
	\item Degrees of freedom are usually chosen by the analyst.  For example, in this work each polygon vertex is chosen to include a displacement degree of freedom in the $x$ and $y$ direction.  Clearly, each element edge has two vertices (points).  A first order polynomial (a line) fits this condition.  Then, values along element edges, for example, can be linearly interpolated.
	\item The number of terms in a polynomial base is equal to $n_k=(k+1)(k+2)$.  For $k=1$, $n_k=6$.  The number of terms in the polynomial base $\mathbf{P}_1$ is 6.
	\item In \eqref{pf2}, rigid body motion is provided by the first three monomials $\mathbf{p}_1, \mathbf{p}_2, \mathbf{p}_3$.	
	\item The number of vertices for a polygon, $n_v$, is equivalent to the number of degrees of freedom, $n_d$, in one spatial direction for a first order ($k=1$) VEM formulation.
\end{enumerate}	

\subsection{Element Stiffness}\label{stiffsect}
Before the stiffness matrix is presented, it is necessary to define a variety of matrices and operators\cite{mengolini,yawVEM}. To begin, recall that the modular matrix for plane stress is
\begin{equation}\label{s2}
	\mathbf{C}=\frac{E_Y}{1-\nu^2}\begin{bmatrix}
		\phantom{0}1 & \phantom{00}\nu & \phantom{00}0\\
		\phantom{0}\nu & \phantom{00}1 & \phantom{00}0\\
		\phantom{0}0 & \phantom{00}0 & \phantom{00}\frac{1-\nu}{2}
		\end{bmatrix},
\end{equation}
and for plane strain is
\begin{equation}\label{s3}
	\mathbf{C}=\frac{E_Y}{(1+\nu)(1-2\nu)}\begin{bmatrix}
		1-\nu & \nu & 0\\
		\nu & 1-\nu & 0\\
		0 & 0 & \frac{1-2\nu}{2}
		\end{bmatrix}\color[rgb]{0,0,0},
\end{equation}
where $E_Y$ is Young's modulus of elasticity, and $\nu$ is Poisson's ratio.

\textbf{Degree of Freedom Operator, $dof_i$.}  Define an operator, $dof_i(\mathbf{p}_j)$, which is found as follows:  first, calculate polynomial vector $\mathbf{p}_j$ evaluated at the vertex coordinates associated with dof $i$, second, take the component of the vector that is directed along dof $i$.  The operator, $dof_i$, extracts the value of its argument, at dof $i$, in the direction of dof $i$, and the dofs (degrees of freedom) range from $1$ to $2n_v$.

\textbf{The $\mathbf{D}$ matrix.}  At coordinates of the degrees of freedom of polygon $E$ polynomial vector values are computed.  These quantities are used to construct matrix $\mathbf{D}$.  The result is
\begin{equation}\label{Dmat}
	\mathbf{D}=\left[ \begin{array}{cccc}
		dof_1(\mathbf{p}_{1}) & dof_1(\mathbf{p}_{2})& \cdots & dof_1(\mathbf{p}_{n_k})\\
		dof_2(\mathbf{p}_{1}) & dof_2(\mathbf{p}_{2})& \cdots & dof_2(\mathbf{p}_{n_k})\\
		\vdots & \vdots & \ddots & \vdots\\
		dof_{2n_v}\mathbf{p}_{1}) & dof_{2n_v}(\mathbf{p}_{2})& \cdots & dof_{2n_v}(\mathbf{p}_{n_k})
		\end{array}\right].
\end{equation}

\textbf{The engineering strain operator, $\bsym{\varepsilon}$.}  An \emph{engineering} strain operator, $\bsym{\varepsilon}$, is applied to the polynomial basis as follows:
\begin{equation}
	\bsym{\varepsilon}(\mathbf{P}_1)=\bsym{\varepsilon}([\mathbf{p}_1 \ \mathbf{p}_1 \ \hdots \ \mathbf{p}_{n_k}])
\end{equation}
Acting on the polynomial base functions, the strain operator produces the matrix,
\begin{equation}\label{ep8}
	\begin{split}
		\bsym{\varepsilon}\left[ \begin{array}{cccc}
		\mathbf{p}_1 & \mathbf{p}_2 & ... & \mathbf{p}_{n_k} 
		\end{array}\right]=\left[ \begin{array}{cccc}
		\partial_{x}p_{1,1} & \partial_{x}p_{2,1} & ... & \partial_{x}p_{n_k,1} \\
		\partial_{y}p_{1,2} & \partial_{y}p_{2,2} & ... & \partial_{y}p_{n_k,2} \\
		\partial_{y}p_{1,1}+\partial_{x}p_{1,2} & \partial_{y}p_{2,1}+\partial_{x}p_{2,2} & ... & \partial_{y}p_{n_k,1}+\partial_{x}p_{n_k,2} 
		\end{array}\right],
	\end{split}
\end{equation}
where in \eqref{ep8}, $p_{i,j}$ represents the $j$th (1st or 2nd) component for polynomial vector $\mathbf{p}_i$.  For example, in 2D the polynomial vectors $\mathbf{p}_i$ are shown in equations \eqref{pf1} and \eqref{pf2}.  \emph{Depending on the context}, the strain operator $\bsym{\varepsilon}$ has a specific use.  For single vectors, like $\mathbf{p}_1$ or $\mathbf{v}^h$, in 2D containing two components, the operator output is $3 \times 1$.  For a group of vectors, like $\mathbf{P}_1$, as shown in \eqref{ep8}, the operator output is $3 \times n_k$, where $n_k=6$ for the first order polynomial base.

\textbf{The engineering stress operator, $\bsym{\sigma}$.}  An \emph{engineering} stress operator, $\bsym{\sigma}(\mathbf{p}_{\alpha})$, is obtained by matrix multiplication
\begin{equation}\label{ne8}
	\begin{split}
		\bsym{\sigma}(\mathbf{p}_{\alpha})&=\mathbf{C}\boldsymbol{\varepsilon}(\mathbf{p}_{\alpha})\\
		\left[ \begin{array}{c}
		\sigma_x(\mathbf{p}_{\alpha})\\
		\sigma_y(\mathbf{p}_{\alpha})\\
		\sigma_{xy}(\mathbf{p}_{\alpha})
		\end{array}\right]&=\mathbf{C}
		\left[ \begin{array}{c}
		\varepsilon_x(\mathbf{p}_{\alpha})\\
		\varepsilon_y(\mathbf{p}_{\alpha})\\
		\gamma_{xy}(\mathbf{p}_{\alpha})
		\end{array}\right].
	\end{split}
\end{equation}
In \eqref{ne8}, the correct modular matrix, $\mathbf{C}$, for plane strain or plane stress is inserted.  The \emph{engineering} stress matrix is formed by using the results of \eqref{ne8}
\begin{equation}\label{ne9}
	\left[\boldsymbol{\sigma}(\mathbf{p}_{\alpha})\right]=\left[ \begin{array}{cc}
		\sigma_x(\mathbf{p}_{\alpha}) & \sigma_{xy}(\mathbf{p}_{\alpha})\\
		\sigma_{xy}(\mathbf{p}_{\alpha}) &\sigma_y(\mathbf{p}_{\alpha})
		\end{array}\right].
\end{equation}
This prepares the way for using \eqref{ne9} in \eqref{ne11}.

\textbf{The $\tilde{\mathbf{B}}$ matrix.}  The $\tilde{\mathbf{B}}$ matrix is formed, for rows $\alpha=1$ to $n_k$ two columns at a time, for $j$ values ranging over vertex numbers $1$ to $n_v$, according to the following formulas:
\begin{equation}\label{ne11}
	\begin{split}
		\tilde{B}_{\alpha (2j-1)}&=\left[ \begin{array}{c}
		1\\
		0
		\end{array}\right]\cdot \left[ \begin{array}{cc}
		\sigma_x(\mathbf{p}_{\alpha}) & \sigma_{xy}(\mathbf{p}_{\alpha})\\
		\sigma_{xy}(\mathbf{p}_{\alpha}) &\sigma_y(\mathbf{p}_{\alpha})
		\end{array}\right]
		\left(\frac{|e_{j-1}|}{2}
		\left[ \begin{array}{c}
		n_{e1}\\
		n_{e2}
		\end{array}\right]_{j-1}+\frac{|e_{j}|}{2}
		\left[ \begin{array}{c}
		n_{e1}\\
		n_{e2}
		\end{array}\right]_j\right)\\ 
		\text{\normalfont and}&\\
		\tilde{B}_{\alpha (2j)}&=\left[ \begin{array}{c}
		0\\
		1
		\end{array}\right]\cdot \left[ \begin{array}{cc}
		\sigma_x(\mathbf{p}_{\alpha}) & \sigma_{xy}(\mathbf{p}_{\alpha})\\
		\sigma_{xy}(\mathbf{p}_{\alpha}) &\sigma_y(\mathbf{p}_{\alpha})
		\end{array}\right]
		\left(\frac{|e_{j-1}|}{2}
		\left[ \begin{array}{c}
		n_{e1}\\
		n_{e2}
		\end{array}\right]_{j-1}+\frac{|e_{j}|}{2}
		\left[ \begin{array}{c}
		n_{e1}\\
		n_{e2}
		\end{array}\right]_j\right).
	\end{split}
\end{equation}
%Typical figure
	\begin{figure}
  \centering
  \epsfig{file=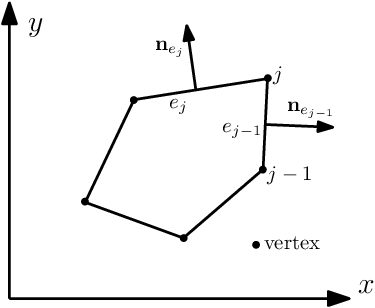,width=.4\textwidth}
  \caption{Single five sided element:  edges, normals, and nodes labeled.}\label{fige}
  \end{figure}
For equation \eqref{ne11}, with reference to Figure \ref{fige}, edge length $j$ is written as $|e_j|$ and the $j$th edge outward unit normal vector components are $[n_{e1} \ n_{e2}]^T_j$.  For vertex $j = 1$, edge length $|e_{j-1}|$ is the length from vertex $j = n_v$ to $j = 1$, where the polygon element has $n_v$ vertices, and normal $n_{e_{j-1}}$ is the outward unit normal to the edge from $j = n_v$ to $j = 1$.  The terms $[1 \ \ 0]^T$ or $[0 \ \ 1]^T$ are dotted with the vector that arises from the right side of \eqref{ne11}, as shown.	

\textbf{The $\breve{\mathbf{B}}$ matrix.}  The $\breve{\mathbf{B}}$ matrix is computed as 
\begin{equation}\label{p16}
	\breve{\mathbf{B}}= \left[ \begin{array}{cccc}
		\breve{\mathbf{b}}_1 & \breve{\mathbf{b}}_2& ... & \breve{\mathbf{b}}_{2n_v} \\
		\end{array}\right],   \quad \quad [3 \times 2n_v]
\end{equation}
where
\begin{equation}\label{p15}
	\breve{\mathbf{b}}_I=\left [ \begin{array}{c}
		\frac{1}{n_v}\sum\limits^{2n_v}_{i=1}\delta_{iI}dof_i(\mathbf{p}_1)\\[0.3cm]
		\frac{1}{n_v}\sum\limits^{2n_v}_{i=1}\delta_{iI}dof_i(\mathbf{p}_2)\\[0.3cm]
		\frac{1}{n_v}\sum\limits^{2n_v}_{i=1}\delta_{iI}dof_i(\mathbf{p}_3)
		\end{array}\right ], \quad \quad \delta_{iI}=\left\{ \begin{array}{c} 
		0 \ \text{for} \ i\neq I\\
		1 \ \text{for} \ i=I.
		\end{array}
		\right.
\end{equation}

\textbf{The $\bar{\mathbf{B}}$ matrix.}  A $\bar{\mathbf{B}}$ matrix is constructed by creating the $\tilde{\mathbf{B}}$ matrix and then replacing the first three rows with the $\breve{\mathbf{B}}$ matrix.  The infinitesimal strain, $\bsym{\varepsilon}(\mathbf{p}_{\alpha})$, is zero for $\alpha=1,2,3$ of the polynomial basis, $\mathbf{P}_1$.  It is for this reason that the beginning three rows of matrix $\tilde{\mathbf{B}}$ given in \eqref{ne11} are populated with $\breve{\mathbf{B}}$ to create the $\bar{\mathbf{B}}$ matrix.  If this was not done the first three rows would be zero (see \cite{mengolini,yawVEM} for more details).  Note, the $\bar{\mathbf{B}}$ matrix described herein should not be confused with the $\bar{\mathbf{B}}$ matrix \cite{Hughes} used with finite elements and incompressibility problems. 

\textbf{The projector matrix, $\tilde{\mathbf{\Pi}}$.}  The projector matrix, $\tilde{\mathbf{\Pi}}$, is constructed as
\begin{equation}\label{eqbbar}
	\tilde{\mathbf{\Pi}}=\mathbf{G}^{-1}\bar{\mathbf{B}},  \quad \quad [n_k \times 2n_v]
\end{equation}
where
\begin{equation}
	\mathbf{G}=\bar{\mathbf{B}}\mathbf{D}.  \quad \quad [n_k \times n_k]
\end{equation}

\textbf{The 2D strain displacement matrix, $\mathbf{B}_2$.}  For 2D VEM the strain displacement matrix is written as
\begin{equation}\label{B2eq}
	\mathbf{B}_2=\bsym{\varepsilon}(\mathbf{P}_1)\tilde{\mathbf{\Pi}}, \quad \quad [3 \times 2n_v]
\end{equation}
\textbf{Engineering strain calculation, $\bsym{\varepsilon}(\mathbf{v}^h)$.}  Using the element vertex displacements, $\mathbf{u}^E$, the engineering strains are computed as shown below:
\begin{equation}\label{ep7}
	\begin{split}
		\bsym{\varepsilon}(\mathbf{v}^h)\approx\bsym{\varepsilon}(\Pi(\mathbf{v}^h))&=\bsym{\varepsilon}\left(\left[ \begin{array}{cccc}
		\mathbf{p}_1 & \mathbf{p}_2 & ... & \mathbf{p}_{n_k} 
		\end{array}\right]\tilde{\bsym{\Pi}}\mathbf{u}^E\right)\\
		&=\bsym{\varepsilon}\left[ \begin{array}{cccc}
		\mathbf{p}_1 & \mathbf{p}_2 & ... & \mathbf{p}_{n_k} 
		\end{array}\right]\tilde{\bsym{\Pi}}\mathbf{u}^E\\
		&=\mathbf{B}_2\mathbf{u}^E. \quad \quad [3 \times 1]
	\end{split}
\end{equation} 
\noindent
\textbf{Remarks}
\begin{enumerate}[(i)]
	\item The strains are ordered according to Voigt notation.  The resulting strain vector $\bsym{\varepsilon}(\mathbf{v}^h)$ contains the two-dimensional engineering strains $\varepsilon_x$, $\varepsilon_y$, $\gamma_{xy}=2\varepsilon_{xy}$.  The size of the strain displacement matrix $\mathbf{B}_2$ in~\eqref{ep7} is $3 \times 2n_v$.
	\item The standard displacement vector ordering is, $\mathbf{u}^E=[u^1_x \ \ u^1_y \ \ u^2_x \ \ u^2_y \cdots u^{n_v}_x \ \ u^{n_v}_y]^T$.
	\item The formulation results in polygon elements with constant strains.
	\item The projection of the virtual element approximation, $\mathbf{v}^h$, onto the polynomial space is implied by the expression $\Pi(\mathbf{v}^h)$.  The expression makes use of the projection operator $\Pi$.
\end{enumerate}

\textbf{Engineering stress calculation, $\bsym{\sigma}(\mathbf{v}^h)$.}  The stress terms needed to calculate the $\tilde{\mathbf{B}}$ matrix are found by using equation \eqref{ne8}.  Yet, actual engineering stresses, calculated during post processing of elasticity problems, are found by using the applicable modular matrix $\mathbf{C}$ and equation \eqref{ep7}. The resulting formula for engineering stresses is
\begin{equation}\label{stress1}
	\bsym{\sigma}(\mathbf{v}^h)\approx \mathbf{C} \bsym{\varepsilon}(\Pi(\mathbf{v}^h))
	=\mathbf{C}\mathbf{B}\mathbf{u}^E.
\end{equation}

\textbf{Element stiffness matrix, $\mathbf{k}_E$.}  Finally, analogous to finite element analysis, the element stiffness matrix for VEM is found found for each polygon element.  The element stiffness for VEM is calculated from two contributions: a consistency part \eqref{kc} and a stability part \eqref{stabmethod2}.
\begin{equation}\label{kE1}
	\mathbf{k}_E=\mathbf{k}_E^c+\mathbf{k}_E^s \quad \quad [2n_v \times 2n_v]
\end{equation}
The consistency part is written as
\begin{equation}\label{kc}
	\mathbf{k}_E^c=tA_E\mathbf{B}_2^T\mathbf{C}\mathbf{B}_2,
\end{equation}
where $\mathbf{C}$ is the 2D plane stress or plane strain elastic modular matrix (see equations \eqref{s2} and \eqref{s3}), $A_E$ is the area of polygon element $E$, and the polygon element has thickness $t$.

The element stability stiffness part, given by Sukumar and Tupek \cite{suku:elastodyn}, is given as
\begin{equation}\label{stabmethod2}
	\mathbf{k}_E^s=(\mathbf{I}-\bsym{\Pi})^T\mathbf{S}^d_E(\mathbf{I}-\bsym{\Pi}).
\end{equation}
where the diagonal matrix $\mathbf{S}^d_E$ is $2n_v \times 2n_v$ in size and is scaled as needed.  The matrix has diagonal terms: $(\mathbf{S}^d_E)_{ii}=\mathsf{max}(\alpha_0 \ \text{\normalfont tr}(\mathbf{C})/m, (\mathbf{k}^c_E)_{ii})$, where $m=3$ in 2D, the modular matrix $\mathbf{C}$ in 2D is the applicable case for plane strain or plane stress, and since elements are constructed from scaled monomials which create diameters on the order of 1, $\alpha_0=1$.

The above stability part of the element stiffness matrix is found to be effective for compressible materials in linear elasticity problems.  It is also appropriate for plasticity problems.  Later, in section \ref{sec:elast}, similar to  Park et al.\cite{park}, a stabilization matrix like \eqref{stabmethod2} is provided for axisymmetric problems of near incompressibility.  

\subsection{Application of External Forces}
Body forces and tractions impose external forces according to equation \eqref{E3}.  Point loads at individual nodes (or vertices) of polygon elements are allowed also.  Standard finite element methods are used to accomplish application of forces.  The following linear form includes the three types of forces mentioned.
\begin{equation}\label{ext1}
	L_E(\mathbf{v}^h)=\int_E \! \mathbf{v}^h\cdot\mathbf{f} \, \mathrm{d}E+\int_{\partial E\cap\Omega_t} \! \mathbf{v}^h\cdot\bar{\mathbf{t}} \, \mathrm{d}\partial E+\sum\limits_{i=1}\mathbf{v}^h(\mathbf{x}_i)\cdot\mathbf{F}_i. 
\end{equation}
Mengolini et al. \cite{mengolini} provides additional discussion on the topic of load application.  In the examples provided later, point loads or pressure loads are used.  An external force vector, $\mathbf{F}_{ext}$, is assembled from the point or pressure loads.  Subsequently, the solution for nodal displacements is found by using the external force vector and the global stiffness.  In anticipation of nonlinear problems, such as plasticity, equilibrium between external and internal forces is accomplished by use of a Newton-Raphson scheme.  

\subsection{Element Internal Forces}
For nonlinear problems, in order to take advantage of an implicit analysis with iterations for equilibrium, local element internal forces are required.  For plastic problems results from the stress integration process are used and a stabilization matrix adjustment is included.  As a result, the internal forces for a single element~$E$ are,
\begin{equation}\label{eif1b}
	\mathbf{q}^E_{int}=A_E \ t \ \mathbf{B}_2^T\bsym{\sigma}+\mathbf{k}_E^s \mathrm{d}\mathbf{u}^E,
\end{equation}
where for element $E$ the polygon area is $A_E$, the engineering stress vector is $\bsym{\sigma}=[\sigma_x \ \sigma_y \ \sigma_{xy}]^T$, and $\mathrm{d}\mathbf{u}^E$ is the latest incremental displacemen.  Note that, for a particular time step in the nonlinear analysis, the stress vector is constant across the given element.
Then, analogous to FEM, the global internal force vector is assembled from the individual element internal force vectors using the assembly operator \cite{Hughes}.  That is,
\begin{equation}\label{eif2}
	\mathbf{F}_{int}=\overset{n_{elem}}{\underset{E=1}{\mathbf{\mathsf{A}}}}\mathbf{q}^{E}_{int}.
\end{equation}

\noindent
\textbf{Remarks}
\begin{enumerate}[(i)]
	\item An implicit nonlinear analysis is possible with the pieces provided above.
	\item The matrix $\mathbf{k}_E$ in \eqref{kE1} is the individual element (tangent) stiffness.  The global (tangent) stiffness matrix, $\mathbf{K}_T$, is assembled from the individual element (tangent) stiffnesses.  It is important to recognize that, although this is similar to FEM, for VEM the number of degrees of freedom for a given element varies depending on the number of vertices in the polygon.  Account of this must be made during assembly.
	\item Equation \eqref{eif2} expresses the global internal force vector, $\mathbf{F}_{int}$.
	\item Like FEM, according to the linear form \eqref{ext1}, the global vector of external forces, $\mathbf{F}_{ext}$, is created.
	\item A residual, $\mathbf{g}=\mathbf{F}_{int}-\lambda \mathbf{F}_{ext}$, for a given load or iteration step is now possible to create when nonlinear problems arise.  The arc length control parameter is $\lambda$, and the residual vector used during equilibrium iterations is given by $\mathbf{g}$.
\end{enumerate}
	
%section 3
\section{Mean Value Coordinates}\label{sec:MVC}
As is seen in a later section, shape function values are needed at the centroid of each VEM element to create axisymmetric virtual elements.  In the VEM formulation, recall that shape function values are not known on the interior of each polygonal element domain.  As a result, a way to calculate the shape functions is needed.  Several options are possible using the nodes of a particular polygon element.  Moving least squares \cite{tb:meshless}, maximum entropy \cite{suku:maxent,suku:maxentreview}, and mean value coordinates (MVC) \cite{Floater} shape functions are all possible solutions.  The least costly, most robust, and easiest to implement are the mean value coordinates shape functions.  They are an effective choice that maintains the benefits of VEM and the ability to handle convex as well as concave polygonal elements.

\begin{figure}
\centering
\epsfig{file=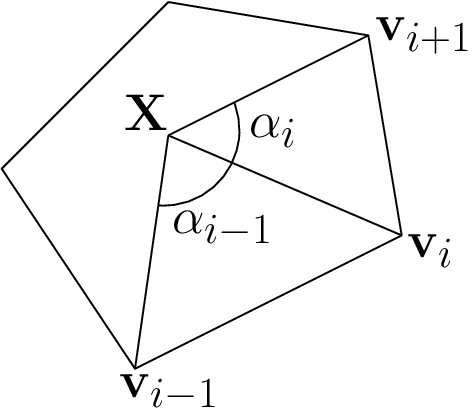,width=.3\textwidth}
\caption{Example Arbitrary Mean Value Coordinates Polygon.}\label{mvcpolygon}
\end{figure}

In Figure~\ref{mvcpolygon} an arbitrary polygon is shown.  For an arbitrary point $\mathbf{X}$ and having the polygon vertices $\mathbf{v}$, the shape function values (and derivatives) are obtained.  As described by Floater \cite{Floater}, shape function values (and derivatives), $\phi_i$, for each node $i$ are calculated as follows:
\begin{equation}
	\phi_i(\mathbf{X})=\frac{w_i(\mathbf{X})}{\sum_{j=1}^n w_j(\mathbf{X})},
\end{equation}
where
\begin{equation}
	w_i(\mathbf{X})=\frac{\tan{(\alpha_{i-1}/2)}+\tan{(\alpha_{i}/2)}}{\|\mathbf{v}_i-\mathbf{X}\|},
\end{equation}
\begin{equation}
	\nabla\phi_i=(\mathbf{R}_i-\sum_{j=1}^n\phi_j\mathbf{R}_j)\phi_i,
\end{equation}
\begin{equation}
	\mathbf{R}_i=\left(\frac{t_{i-1}}{t_{i-1}+t_i}\right)\frac{\mathbf{c}_{i-1}^{\perp}}{\sin{\alpha_{i-1}}}+
	\left(\frac{t_{i-1}}{t_{i-1}+t_i}\right)\frac{\mathbf{c}_{i}^{\perp}}{\sin{\alpha_{i}}}+\frac{\mathbf{e_i}}{r_i}.
\end{equation}
The following definitions apply for the above mean value coordinates formulas:
\begin{align}
	t_i&=\tan{(\alpha_i/2)}, \quad 0<\alpha_i<\pi\\
	r_i&=||\mathbf{v}_i-\mathbf{X}||\\
	\mathbf{e}_i&=\frac{(\mathbf{v}_i-\mathbf{X})}{||\mathbf{v}_i-\mathbf{X}||}\\
	\mathbf{c}_i&=\frac{\mathbf{e}_i}{r_i}-\frac{\mathbf{e}_{i+1}}{r_{i+1}}\\
	\text{for a vector}& \ \mathbf{a}=(a_1,a_2)\in \Re^2, \ \text{let} \ \mathbf{a}^{\perp}=(-a_2,a_1).
\end{align}

%section 4 The axisymmetric formulation
\section{The axisymmetric formulation}\label{sec:axi}
With the results of sections \ref{sec:reviewVEM} and \ref{sec:MVC}, the pieces necessary to construct axisymmetric VEM are in hand.  A non-isoparametric formulation is presented since VEM and MVC shape functions are in global coordinates.  For the following derivation, the vector of shape functions is denoted as $\bsym{\phi}=[\phi_1 \ \phi_2 \ \cdots \ \phi_{n_v}]^T$.

\begin{figure}
\centering
\epsfig{file=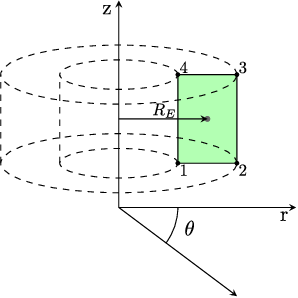,width=.4\textwidth, clip=} %loaded nice pdf into inkscape and converted it to eps there, then use clip=
\caption{Axisymmetric element.}\label{axifig}
\end{figure}

\subsection{Axisymmetric displacements}%approximation u
Consider Figure \ref{axifig} with an axisymmetric element shown.  The example element shown in the figure has an axisymmetric cross-section with 4 nodes.  In general, the cross-section has an arbitrary polygonal cross-section with nodes $1$ to $n_v$.  For the cross-section the nodal displacements in the radial direction and $z$ direction are represented in terms of nodal shape functions as follows:
\begin{equation}
	\mathbf{u}=\left\{\begin{array}{c}		
		u_r \\
		u_z
		\end{array}\right\}=\left\{\begin{array}{c}		
		\bsym{\phi}^T\mathbf{d}_r \\
		\bsym{\phi}^T\mathbf{d}_z
		\end{array}\right\},
\end{equation}
where $\mathbf{d}_r=[u_r^1 \ u_r^2 \ \cdots \ u_r^{n_v}]^T$ and $\mathbf{d}_z=[u_z^1 \ u_z^2 \ \cdots \ u_z^{n_v}]^T$.

\subsection{Axisymmetric cross-sectional strains}%In plane strains
The 2D coordinates are relabeled in terms of the axisymmetric variables $r\equiv x$, $z \equiv y$.  As a result, for a given element, the axisymmetric cross-sectional strains are expressed as
\begin{equation}
	\bsym{\varepsilon}=\left\{\begin{array}{c}		
		\varepsilon_r \\
		\varepsilon_z \\
		\gamma_{rz}
		\end{array}\right\}=\left\{\begin{array}{c}		
		\frac{\partial u_r}{\partial r} \\
		\frac{\partial u_z}{\partial z} \\
		\frac{\partial u_z}{\partial r}+\frac{\partial u_r}{\partial z}
		\end{array}\right\}=
		\setlength\arraycolsep{4pt}
		\left[ \begin{array}{ccccccc}		
		\frac{\partial \phi_1}{\partial r} & 0 & \frac{\partial \phi_2}{\partial r} & 0 & ... & \frac{\partial \phi_{n_v}}{\partial r} & 0 \\
		0 & \frac{\partial \phi_1}{\partial z} & 0 & \frac{\partial \phi_2}{\partial z} &... & 0 & \frac{\partial \phi_{n_v}}{\partial z} \\
		\frac{\partial \phi_1}{\partial z} &\frac{\partial \phi_1}{\partial r} & \frac{\partial \phi_2}{\partial z} &\frac{\partial \phi_2}{\partial r} & ... & \frac{\partial \phi_{n_v}}{\partial z} &\frac{\partial \phi_{n_v}}{\partial r}		
		\end{array}\right]
		\left\{\begin{array}{c}		
		u_r^1 \\
		u_z^1 \\
		\vdots\\
		u_r^{n_v} \\
		u_z^{n_v} \\
		\end{array}\right\}=\mathbf{B}_2\mathbf{d},
\end{equation}
where the vector $\mathbf{d}$ organizes $r$ and $z$ displacements as shown.  Notice the strain displacement matrix, $\mathbf{B}_2$.  Aside from coordinate direction relabeling, it is conceptually identical to the 2D elasticity VEM strain displacement matrix determined in section \ref{stiffsect}, equation \eqref{B2eq}.  The shape functions themselves are not known, but the strain displacement matrix is computable from the VEM formulation.

\subsection{The tangential or circumferential strain}%The tangential or circumferential strain
For a given element in the VEM discretization, the tangential or circumferential strain is
\begin{equation}
	\begin{split}
		\varepsilon_t &=\frac{C_f-C_o}{C_o}=\frac{2\pi(R_E+u_r)-2\pi R_E}{2\pi R_E}=\frac{u_r}{R_E}\\
			&=\frac{1}{R_E}\bsym{\phi}^T\mathbf{d}_r=\frac{1}{R_E}\left[\phi_1 \ 0 \ \cdots  \phi_{n_v} \ 0 \right]\mathbf{d} \\
			&=\mathbf{B}_t\mathbf{d},
	\end{split}
\end{equation}
where $C_o$ is the original circumference, $C_f$ is the final (or current) circumference, and $R_E$ is the distance from the axis of symmetry to the element centroid (see Figure~\ref{axifig}), and the shape functions associated with tangential strain are evaluated at the element centroid coordinates.  Note also that the strain displacement matrix $\mathbf{B}_t$ requires shape function values.  The shape function value are not available from the VEM formulation.  Hence, an alternative means of calculating shape function values is needed.  The mean value coordinates (MVC) shape functions of section \ref{sec:MVC} are used.  

\subsection{The axisymmetric strain displacement matrix}
The strain displacement matrix for the axisymmetric case is constructed by combining the results from the in-plane case, $\mathbf{B}_2$, which comes from the VEM formulation, and the circumferential or tangential strain case, $\mathbf{B}_t$, which is constructed by using MVC shape functions.  The result is the axisymmetric strain displacement matrix
\begin{equation}
		\mathbf{B}=\left[\begin{array}{c}
			\mathbf{B}_2\\
			\mathbf{B}_t
		\end{array}
		\right]=
		\setlength\arraycolsep{4pt}
		\left[ \begin{array}{ccccccc}		
		\frac{\partial \phi_1}{\partial r} & 0 & \frac{\partial \phi_2}{\partial r} & 0 & ... & \frac{\partial \phi_{n_v}}{\partial r} & 0 \\
		0 & \frac{\partial \phi_1}{\partial z} & 0 & \frac{\partial \phi_2}{\partial z} &... & 0 & \frac{\partial \phi_{n_v}}{\partial z} \\
		\frac{\partial \phi_1}{\partial z} &\frac{\partial \phi_1}{\partial r} & \frac{\partial \phi_2}{\partial z} &\frac{\partial \phi_2}{\partial r} & ... & \frac{\partial \phi_{n_v}}{\partial z} &\frac{\partial \phi_{n_v}}{\partial r} \\
		\frac{\phi_1}{R_E} & 0 & \frac{\phi_2}{R_E} & 0 & ... & \frac{\phi_{n_v}}{R_E} & 0
		\end{array}\right].
\end{equation}

%%%%Axisymmetric C
\subsection{Axisymmetric modular matrix and stresses}
The axisymmetric modular matrix is
\begin{equation}\label{s3b}
	\mathbf{C}=\zeta\begin{bmatrix}
		1 & \frac{\nu}{1-\nu} & 0 & \frac{\nu}{1-\nu}\\
		\frac{\nu}{1-\nu} & 1 & 0 & \frac{\nu}{1-\nu}\\
		0 & 0 & \frac{1-2\nu}{2(1-\nu)} & 0\\
		\frac{\nu}{1-\nu} & \frac{\nu}{1-\nu} & 0 & 1
		\end{bmatrix}\color[rgb]{0,0,0},
\end{equation}\\
where $\zeta=\frac{E_Y(1-\nu)}{(1+\nu)(1-2\nu)}$.

The axisymmetric stresses are expressed as
\begin{equation}\label{ne8b}
	\begin{split}
		\bsym{\sigma}&=\mathbf{C}\boldsymbol{\varepsilon}\\
		\left[ \begin{array}{c}
		\sigma_r\\
		\sigma_z\\
		\sigma_{rz}\\
		\sigma_t
		\end{array}\right]&=\mathbf{C}
		\left[ \begin{array}{c}
		\varepsilon_r\\
		\varepsilon_z\\
		\gamma_{rz}\\
		\varepsilon_t
		\end{array}\right].
	\end{split}
\end{equation}
%%%%  

% Section 5
\section{Axisymmetric VEM - Elasticity}\label{sec:elast}
For problems of axisymmetric elasticity, using VEM notation, the consistent part of the element stiffness is expressed as
\begin{equation}
	\mathbf{k}_E^c=2 \pi R_E A_E\mathbf{B}^T \mathbf{C} \mathbf{B},
\end{equation}
where the 2D element volume of $V_E=A_Et$ is replaced with the axisymmetric volume, $V_E=2\pi R_E A_E$.

The stability part is expressed like the 2D elasticity case as
\begin{equation}\label{elastaxistab}
	\mathbf{k}_E^s=(\mathbf{I}-\bsym{\Pi})^T\mathbf{S}^d_E(\mathbf{I}-\bsym{\Pi}).
\end{equation}
However, with an eye toward compressible \emph{and} near incompressible problems, a procedure similar to \cite{park}, is used.  Consequently, the diagonal matrix $\mathbf{S}^d_E$ is defined to have diagonal terms: $(\mathbf{S}^d_E)_{ii}=(\mathbf{k}^c_{E,\mu})_{ii}$, where $\mathbf{k}^c_{E,\mu}$ is the consistency matrix constructed with the $\mu=\frac{E_Y}{2(1+\nu)}$ part \cite{Hughes} of $\mathbf{C}$.  That is
\begin{equation}\label{Cmu}
	\mathbf{C}_{\mu}=
	\setlength\arraycolsep{4pt}
	\mu\begin{bmatrix}
		2 & 0 & 0 & 0\\
		0 & 2 & 0 & 0\\
		0 & 0 & 1 & 0\\
		0 & 0 & 0 & 2
		\end{bmatrix}\color[rgb]{0,0,0}.
\end{equation}\\

The element stiffness matrix is then formed in the usual way as
\begin{equation}
	\mathbf{k}_E=\mathbf{k}_E^c+\mathbf{k}_E^s.
\end{equation}

Special care is required to properly apply loads in axisymmetric problems.  For an axisymmetric load at node $i$ the applied force is
\begin{equation}
	F_j=2\pi r_j w_c,
\end{equation}
where $r_j$ is the radial distance to node $j$ from the axis of symmetry, and $w_c$ is the load per unit distance along the circumference associated with $r_j$.  For an axisymmetric pressure at node $j$ the applied force is
\begin{equation}
	F_j=2\pi r_j h p,
\end{equation}
where $p$ is the pressure, and $h$ is the distance along the surface tributary to node $j$.  For example, if the surface is parallel to the $z$ axis then $h=|z_{j+1}-z_{j-1}|/2$, where $z_{j+1}$ is the $z$ coordinate of the node above $z_j$ and $z_{j-1}$ is the $z$ coordinate of the node below $z_j$.

%Section 6
\section{Axisymmetric VEM - Plasticity}\label{sec:plast}
Material nonlinearities like plasticity are easily included.  For such a case, material properties must be updated during each load step.  The local axisymmetric stiffness is the same as the case for elasticity except that a consistent elasto-plastic modular matrix, $\mathbf{C}_{ep}$, is inserted into the expressions for $\mathbf{k}_E^c$ and $\mathbf{k}_E^s$.  Then the matrix $\mathbf{k}_{E}$ becomes
\begin{equation}\label{plastaxistab}
	\mathbf{k}_E= \mathbf{k}_E^c+\mathbf{k}_E^s=2\pi R_E A_E\mathbf{B}^T \mathbf{C}_{ep} \mathbf{B}+\mathbf{k}_E^s(\mathbf{C}_{ep}).
\end{equation}
As strains evolve during each load step, stresses and $\mathbf{C}_{ep}$ are updated according to the 3D $J2$ plasticity formulation with radial return (see Simo and Taylor \cite{Simo3}, and Simo and Hughes \cite{Simo2}).  The stabilization matrix here is a function of $\mathbf{C}_{ep}$ and is formed by taking the diagonal terms of the elasto-plastic consistency matrix, $\mathbf{k}_E^c$, according to the procedure for stabilization matrix \eqref{stabmethod2} so that $(\mathbf{S}^d_E)_{ii}=(\mathbf{k}^c_{E})_{ii}$.  All other formulas remain the same.

The internal force vector for an individual element is
\begin{equation}\label{eifplast}
	\mathbf{q}^E_{int}=2\pi R_E \ A_E \ \mathbf{B}^T\bsym{\sigma}+\mathbf{k}_E^s \mathrm{d}\mathbf{u}^E,
\end{equation}
where $\bsym{\sigma}=[\sigma_{rr} \ \sigma_{zz} \ \tau_{rz} \ \sigma_{tt} ]^T$ is the vector of element stresses obtained from the plasticity stress integration algorithm, and $\mathrm{d}\mathbf{u}^E$ is the latest incremental displacement.

% Section 7
\section{Numerical Implementation}\label{sec:numimplementation}
For the nonlinear problems involving plasticity, implicit Newton-Raphson iterations are used with radial return at the constitutive level for each element.  Global equilibrium is enforced by an arc-length path following scheme \cite{Crisfieldv1}.  As a result, the nonlinear analysis is carried out with typical ingredients: (i) a path following scheme (arc-length method), (ii) global external force vector, (iii) global internal force vector, and (iv) consistent global tangent stiffness matrix.  For all example elasticity and plasticity problems the form of stabilization used is described by equations \eqref{elastaxistab} and \eqref{plastaxistab}, respectively.  All polygonal mesh generation is accomplished by using Polymesher \cite{talischi}.

% Section 8
\section{Numerical Results}\label{sec:numresults}
A variety of simulations are provided to illustrate the utility of the axisymmetric VEM.  In particular, representative results of static elastic and plastic simulations, using consistent units, are provided.  Results are compared to theoretical solutions or benchmark finite element solutions.  For contour plots of stresses the calculated element stresses are used directly.  For plots of stress versus radial distance, nodal stresses are shown and are calculated based on the average of element stresses common to a given node.  

Convergence rates of axisymmetric VEM are computed using the $L_2$ displacement discrete error measure \cite{Chen3}.  The $L_2$ displacement error measure is calculated as
\begin{equation}
	||\mathbf{u}-\mathbf{u}^h||_{L_2(\Omega)}=\sqrt{\sum\limits_E \int\limits_E |\mathbf{u}-\mathbf{u}^h|^2 \,\mathrm{d}\mathbf{x}}.
\end{equation}
Consequently, the \emph{relative} $L_2$ displacement error is calculated as
\begin{equation}
	\frac{||\mathbf{u}-\mathbf{u}^h||_{L_2(\Omega)}}{||\mathbf{u}||_{L_2(\Omega)}}.
\end{equation}

\subsection{Cylinder Under Internal Pressure - Linear Elastic Lam\'{e} Problem, Case of Plane Strain}
\begin{figure}
\centering
\mbox{
\subfigure[]{\epsfig{file=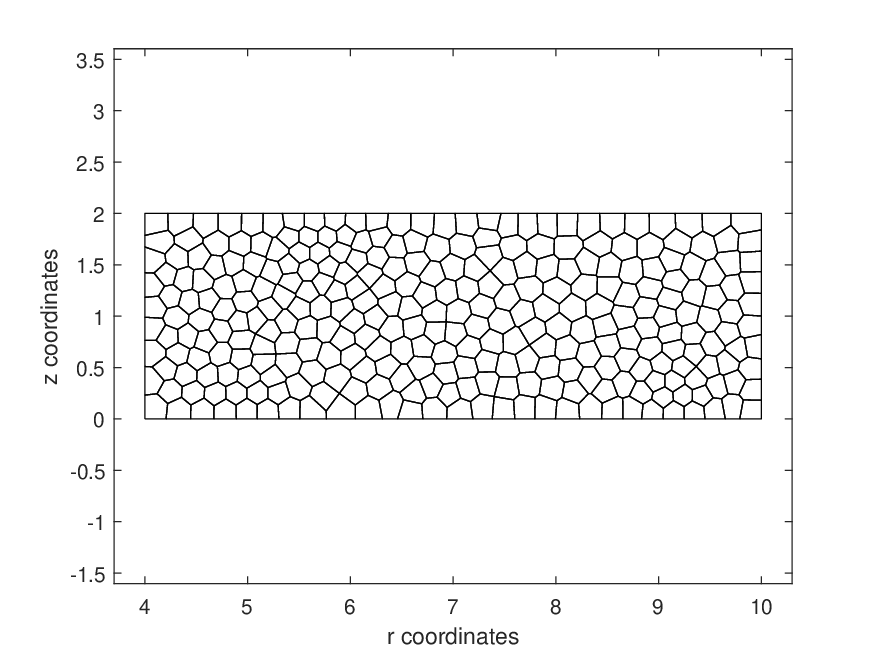,width=0.48\textwidth}}
\subfigure[]{\epsfig{file=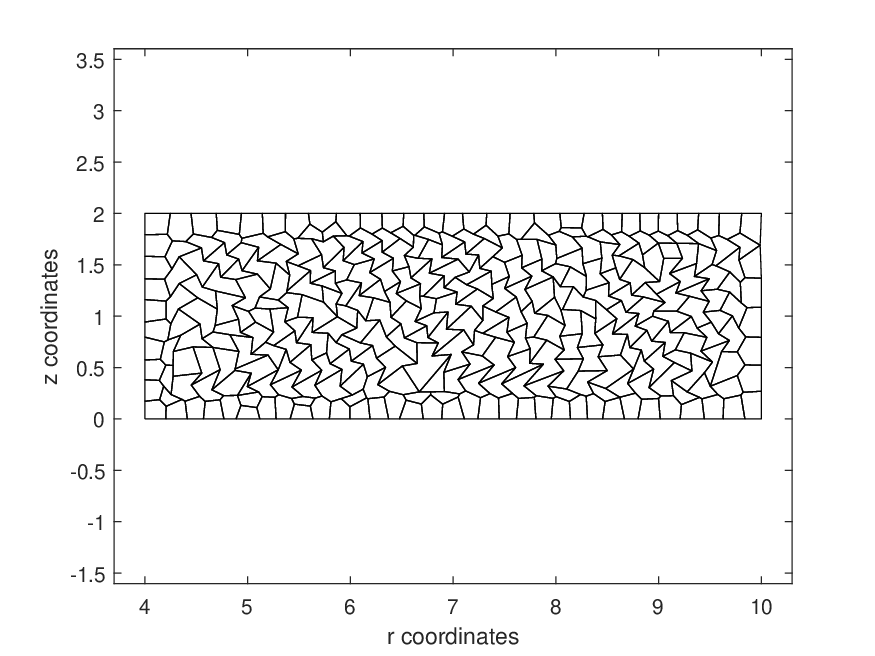,width=0.48\textwidth}}
}

%\vspace*{-0.2in}
\mbox{
\subfigure[]{\epsfig{file=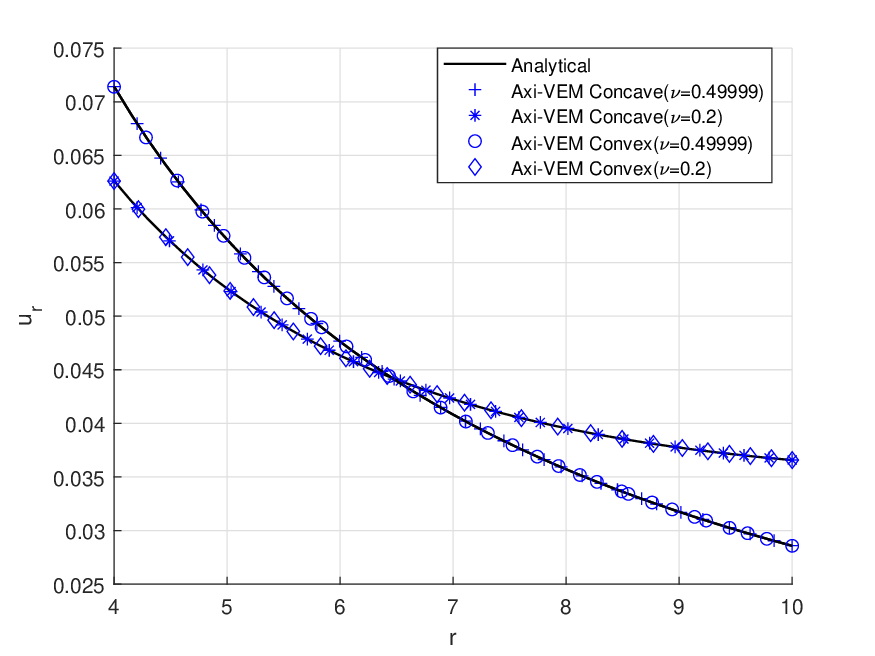,width=0.48\textwidth}}
\subfigure[]{\epsfig{file=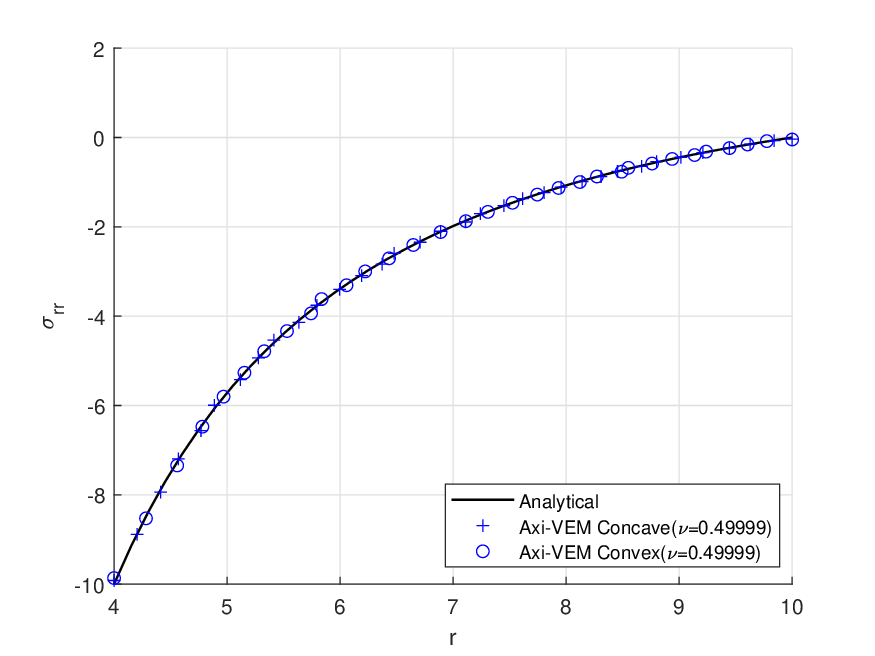,width=0.48\textwidth}}
}

%\vspace*{-0.2in}
\mbox{
\subfigure[]{\epsfig{file=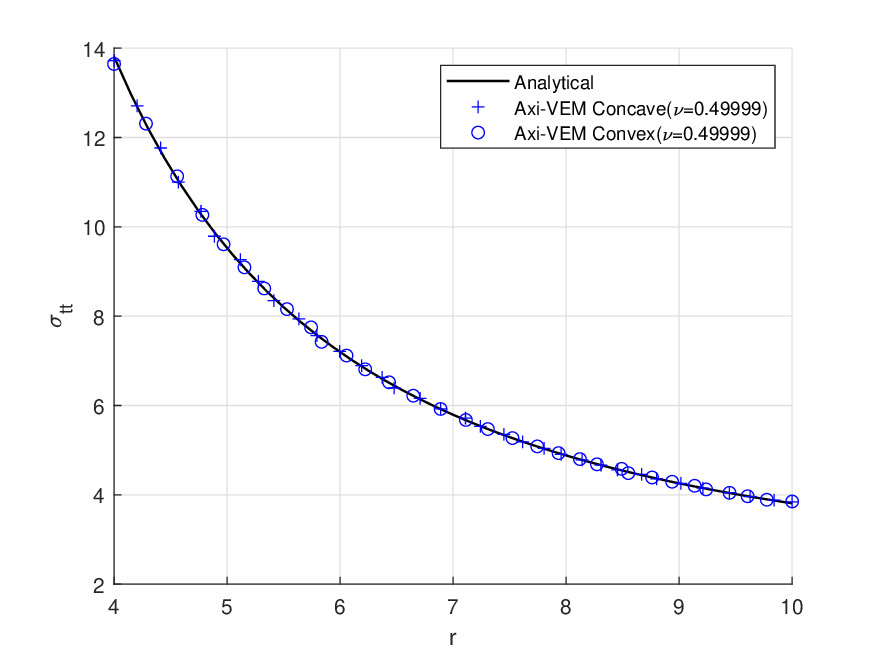,width=0.48\textwidth}}
\subfigure[]{\epsfig{file=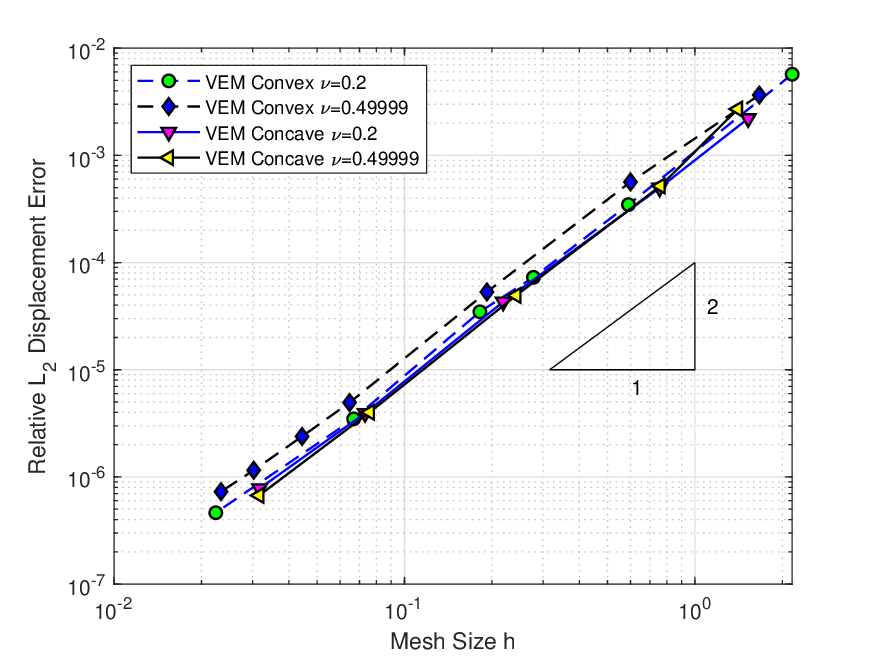,width=0.48\textwidth}}
}

\vspace*{-0.1in}
\caption{Plain Strain Elastic Cylinder Under Internal Pressure:  (a)~VEM discretization (300 convex elements); (b)~VEM discretization (300 mostly concave elements); (c)~Radial displacement versus radial distance; (d)~Radial stress, $\sigma_{rr}$, versus radial distance; and (e)~Tangential stress, $\sigma_{tt}$, versus radial distance, (f)~Relative $L_2$ displacement error versus mesh size $h$.}\label{numericallinearcyl}
\end{figure}

For the condition of plane strain a linear elastic cylinder is subject to internal pressure, $p=10$.  The cylinder has inner radius $a=4$, and outer radius, $b=10$.  The modulus of elasticity $E_Y=1000$ and Poisson's ratio $\nu=0.2$ or $\nu=0.49999$.  Results are presented for the case of plane strain in the $r-\theta$ plane.  Boundary conditions are set so that nodes at the top and bottom of the discretization are restrained in the $z$ direction causing $\varepsilon_z=0$ everywhere. Plots for Figures \ref{numericallinearcyl}abcde, are accomplished with 300 VEM elements.  The theoretical solution \cite{Timoshenko}, for radial displacements and stresses, is
\begin{equation}
	\begin{split}
		u_r&=\frac{pa^2(1+\nu)(b^2+r^2(1-2\nu))}{rE_Y(b^2-a^2)}\\
		\sigma_r&=\frac{pa^2}{b^2-a^2}\left(1-\frac{b^2}{r^2}\right)\\
		\sigma_z&=\frac{2\nu pa^2}{b^2-a^2}\\
		\sigma_{rz}&=0\\
		\sigma_t&=\frac{pa^2}{b^2-a^2}\left(1+\frac{b^2}{r^2}\right).
		\end{split}
\end{equation}
For problems of near incompressibility the results for Possion's ratio equal to 0.49999 are handled without trouble.  Convex and concave results are accomplished with indistinguishable results.  The convergence characteristics shown in Figure~\ref{numericallinearcyl}f are within expected ranges for both compressible and near incompressible values of Poisson's ratio and are not affected by concave or convex elements.

\subsection{Sphere Under Internal Pressure - Linear Elastic Lam\'{e} Problem}
\begin{figure}
\centering
\mbox{
\subfigure[]{\epsfig{file=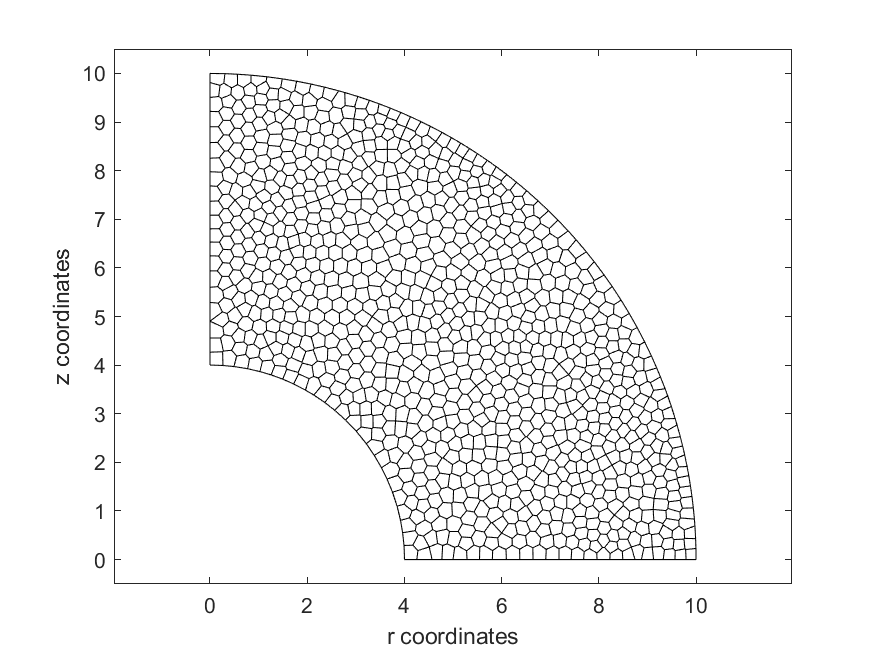,width=0.48\textwidth}}
\subfigure[]{\epsfig{file=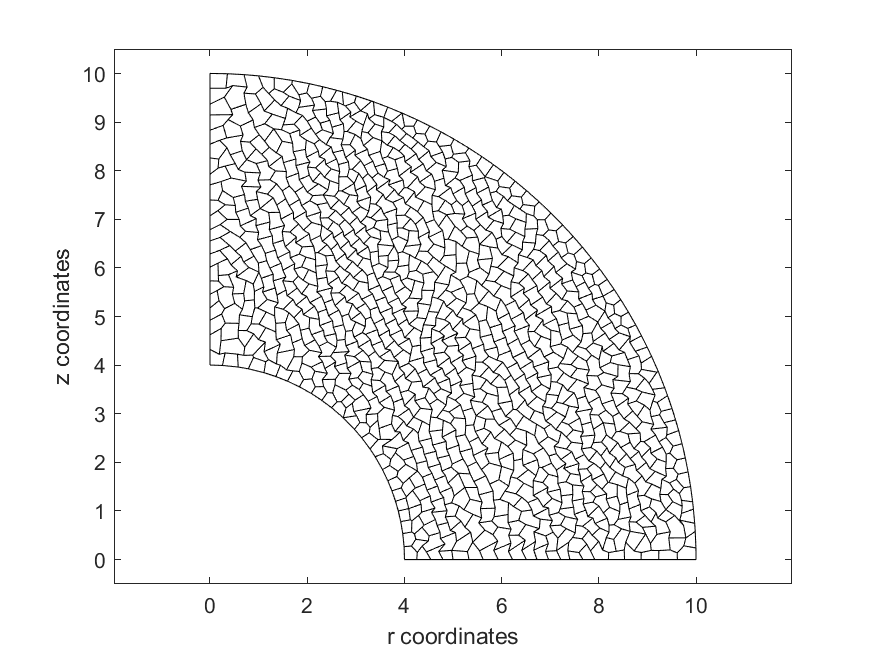,width=0.48\textwidth}}
}

%\vspace*{-0.2in}
\mbox{
\subfigure[]{\epsfig{file=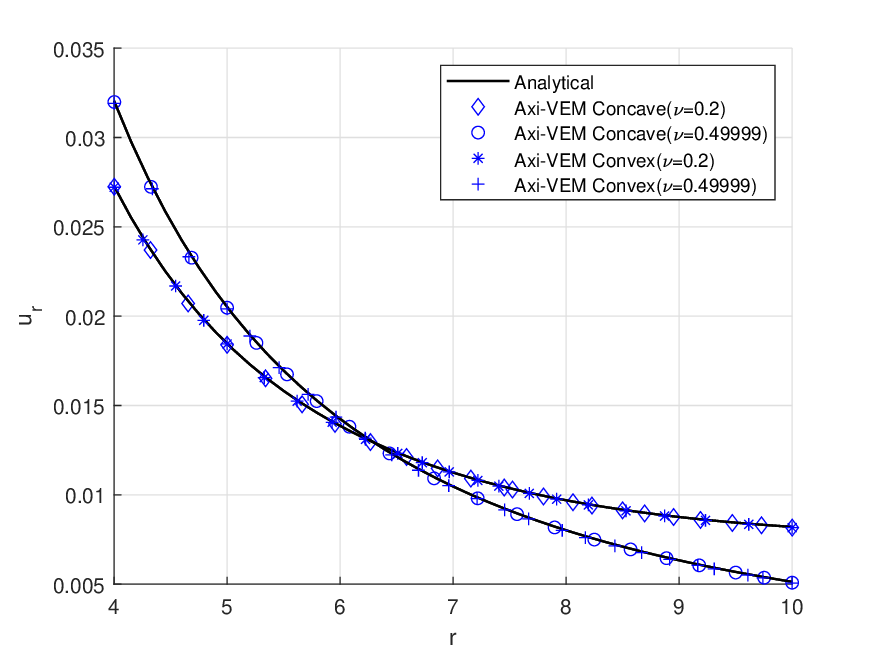,width=0.48\textwidth}}
\subfigure[]{\epsfig{file=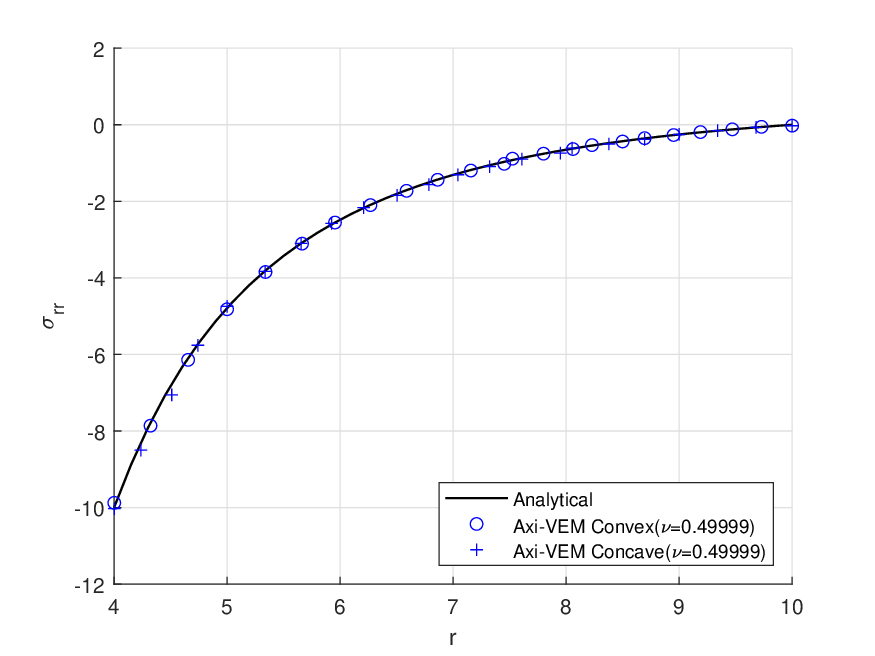,width=0.48\textwidth}}
}

%\vspace*{-0.2in}
\mbox{
\subfigure[]{\epsfig{file=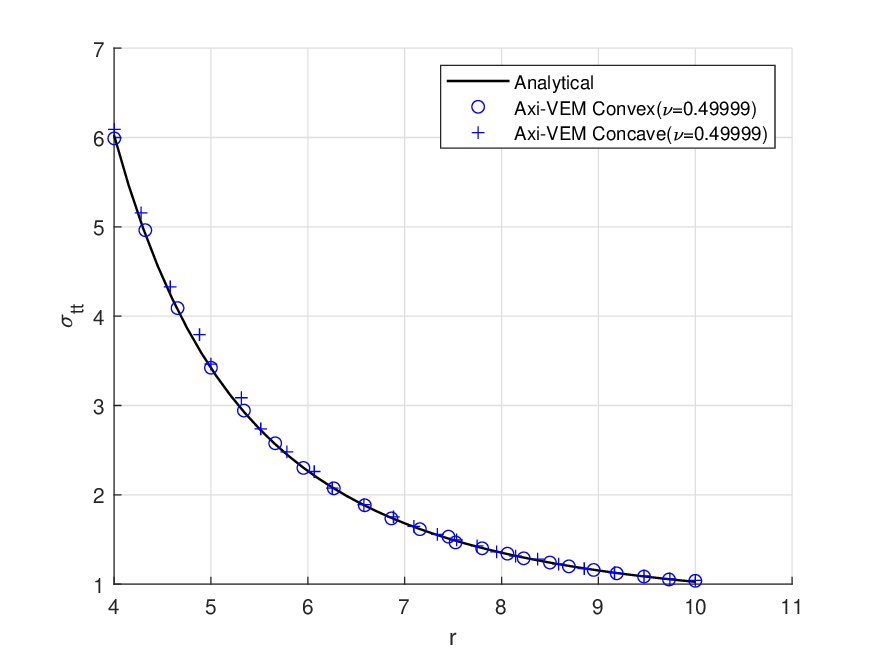,width=0.48\textwidth}}
\subfigure[]{\epsfig{file=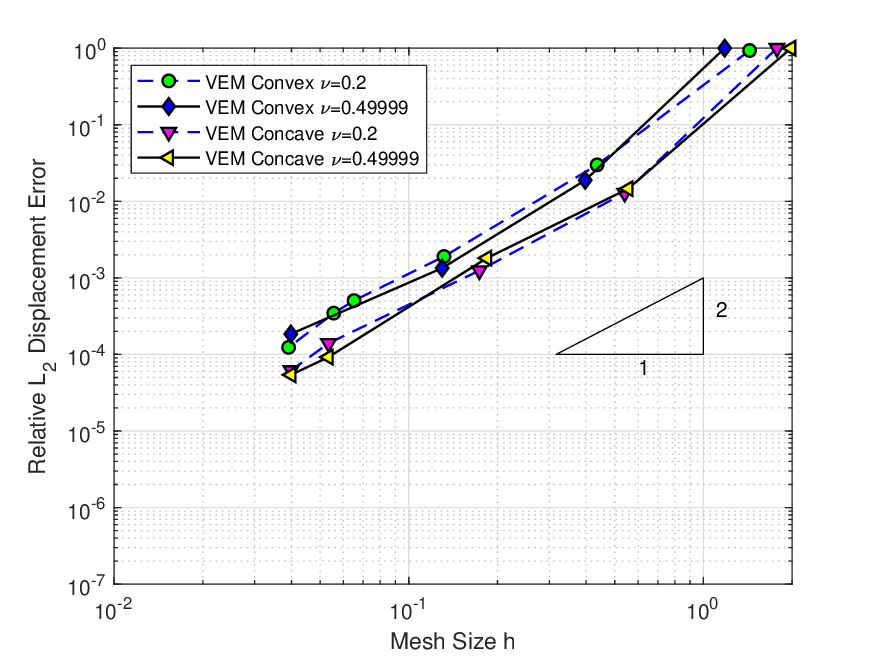,width=0.48\textwidth}}
}

\vspace*{-0.1in}
\caption{Linear Elastic Sphere Under Internal Pressure:  (a)~VEM discretization (1000 convex elements); (b)~VEM discretization (1000 mostly concave elements); (c)~Radial displacement versus radial distance; (d)~Radial stress, $\sigma_{rr}$, versus radial distance; and (e)~Tangential stress, $\sigma_{tt}$, versus radial distance, (f)~Relative $L_2$ displacement error versus mesh size $h$.}\label{numericallinearsphere}
\end{figure}

In Figures~\ref{numericallinearsphere}ab, the top half of an axisymmetric sphere is modeled with 1000 polygonal elements.  The applied internal pressure is $p=10$.  The sphere internal radius $a=4$, the external radius $b=10$, modulus of elasticity $E_Y=1000$, and Poisson's ratio $\nu=0.2$ or $\nu=0.49999$.  Vertical supports are provided at all nodes for which $z=0$, and horizontal supports are provided at all nodes for which $r=0$.  The theoretical solution \cite{Timoshenko}, for radial displacements and stresses, is
\begin{equation}
	\begin{split}
		u_r&=\frac{pa^3}{Er^2(b^3-a^3)}\left[\frac{(1-\nu)(2r^3+b^3)}{2}+\nu(b^3-r^3)\right]\\
		\sigma_r&=\frac{pa^3(b^3-r^3)}{r^3(a^3-b^3)}\\
		\sigma_{rz}&=0\\
		\sigma_t&=\frac{-pa^3(2r^3+b^3)}{2r^3(a^3-b^3)}.\\
		\end{split}
\end{equation}

\subsection{Elastic Circular Plate with Center Point Load}
\begin{figure}
\centering
\mbox{
\subfigure[]{\epsfig{file=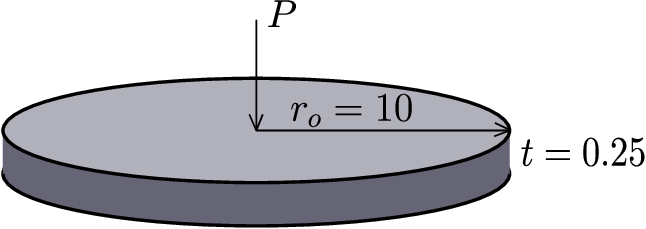,width=0.40\textwidth,clip=}} \hspace{0.3in}
\subfigure[]{\epsfig{file=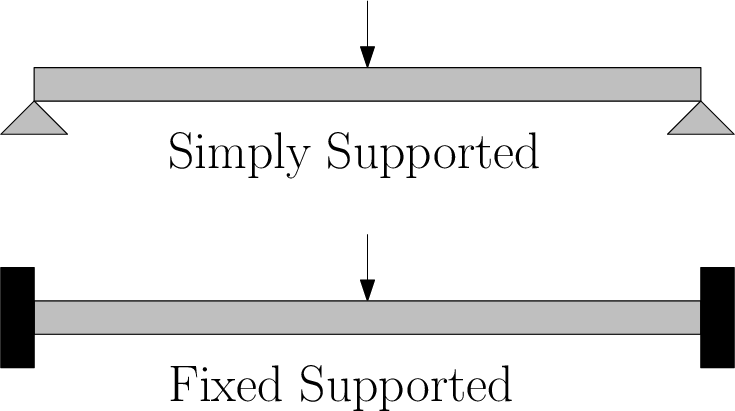,width=0.35\textwidth}}
}

%\vspace*{-0.2in}
\mbox{
\subfigure[]{\epsfig{file=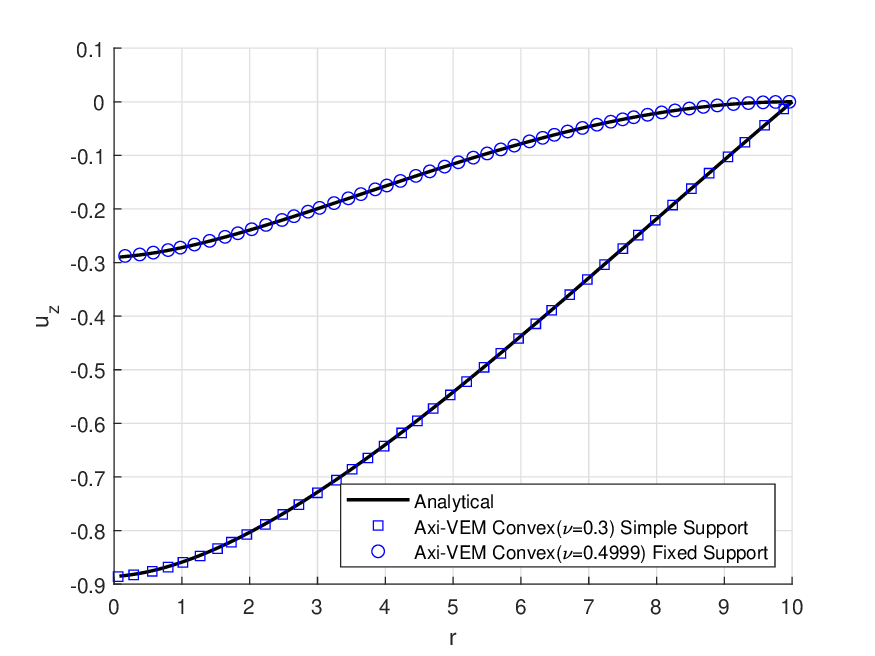,width=0.48\textwidth}}
\subfigure[]{\epsfig{file=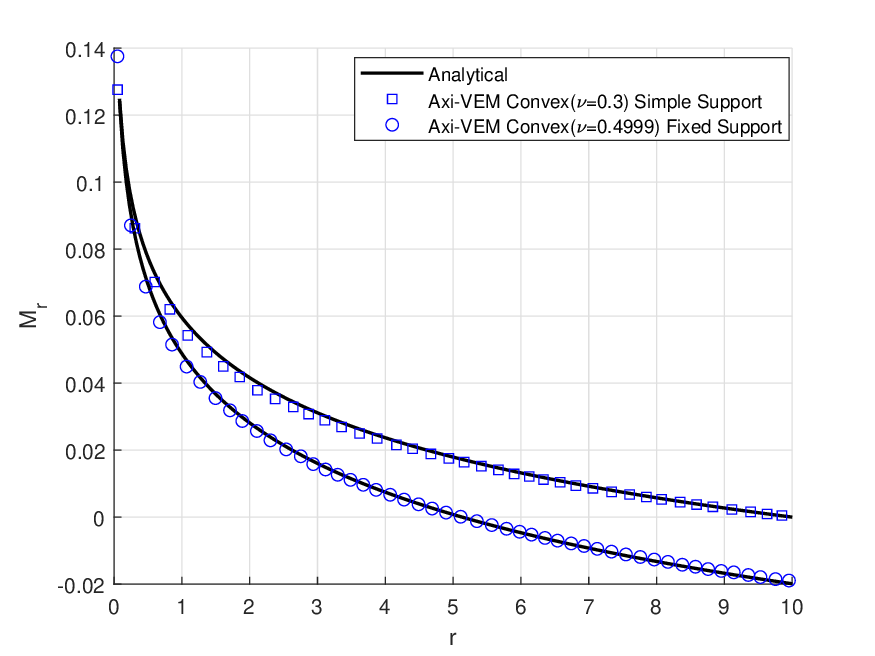,width=0.48\textwidth}}
}

%\vspace*{-0.2in}
\mbox{
\subfigure[]{\epsfig{file=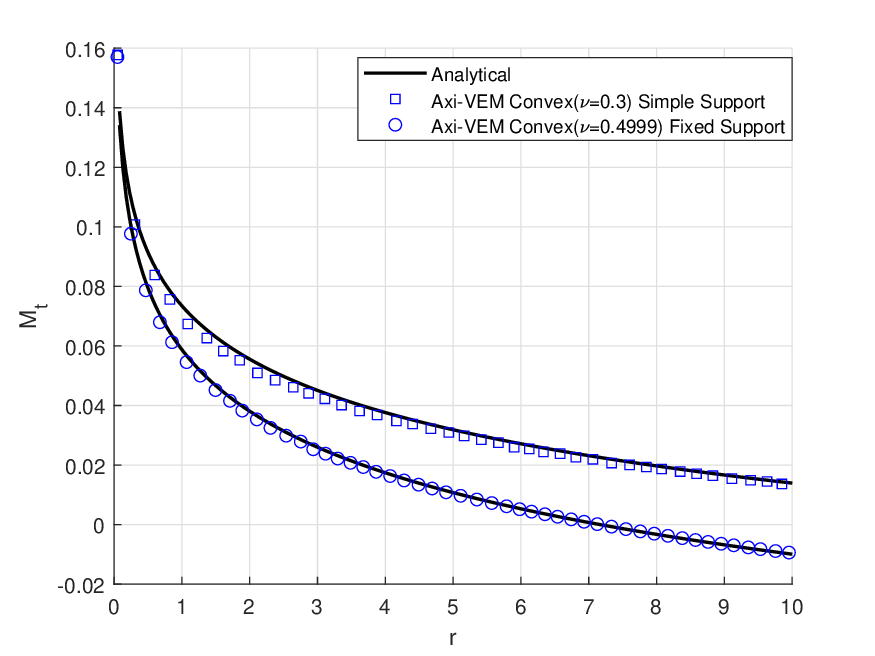,width=0.48\textwidth}}
%\subfigure[]{\epsfig{file=./figs/sigmasycant.eps,width=0.48\textwidth}}
}

\vspace*{-0.1in}
\caption{Circular plate with point load ($P=0.25$) at center:  (a)~Circular plate geometry and loading; (b)~Simple support and fixed support boundary condition cases; (c)~Center displacement, $u_z$, versus radial coordinate, $r$; (d)~Radial moment, $M_r$, versus radial coordinate, $r$; (e)~Tangential moment, $M_t$, versus radial coordinate, $r$.}
\label{numericalplateptload}
\end{figure}

An elastic circular plate is loaded with a point load at the center, as shown in Figure~\ref{numericalplateptload}a.  The radius of the plate is $
r_o=10$.  The plate has thickness $t=0.25$ and modulus of elasticity $E_Y=1000$.  Poisson's ratio used in the various results are indicated in Figures~\ref{numericalplateptload}cde.  In particular, for some results Poisson's ratios as high as 0.49999 are used to demonstrate that the formulation is not susceptible to locking for near incompressible problems.  Results for displacements, radial moments, and tangential moments versus radial distance, $r$, are shown to be in good agreement with analytical solutions.  Furthermore, Figures~\ref{numericalplateptload}cde illustrate results for plates with simply supported and fixed supported conditions.  These problems are accomplished with an axisymmetric model using 18000 convex polygonal elements.  

A study of normalized displacements and various Poisson ratios is provided in Table~\ref{numplateptloadTab}.  For all results of the table, the circular plate is modeled with 9000 convex or mostly concave elements, as indicated.  The results are in good agreement with the theoretical displacements which account for both bending and shear displacements.  The axisymmetric VEM results have error of at most 3.76 \% and for most cases less than 2\%.  The axisymmetric VEM formulation using concave or convex elements provides nearly indistinguishable results.

For the circular plate loaded by a point load at the center, the center displacement (bending plus shear) is $u_z$, the radial moment is $M_r$, and the tangential moment is $M_t$.  The theoretical solutions \cite{Timoshenko2,Ugural} for an elastic circular plate with edges simply supported are
\begin{equation}
	\begin{split}
		u_z&=\frac{P}{16\pi D}\left[\frac{3+\nu}{1+\nu}(r_o^2-r^2)+2r^2\ln{\frac{r}{r_o}}\right]-\frac{Pt^2}{8\pi D(1-\nu)}\ln{\frac{r}{r_o}}\\
    M_r&=\frac{P}{4\pi}(1+\nu)\ln{\frac{r_o}{r}}\\
    M_t&=\frac{P}{4\pi}\left[(1+\nu)\ln{\frac{r_o}{r}}+1-\nu\right]\\
		D&=\frac{E_Yt^3}{12(1-\nu^2)}
	\end{split}
\end{equation}
For the case of an elastic circular plate with edges fixed supported, the theoretical solutions \cite{Timoshenko2,Ugural} are
\begin{equation}
	\begin{split}
		u_z&=\frac{Pr^2}{8\pi D}\ln{\frac{r}{r_o}}+\frac{P(r_o^2-r^2)}{16\pi D}-\frac{Pt^2}{8\pi D(1-\nu)}\ln{\frac{r}{r_o}}\\
    M_r&=\frac{P}{4\pi}\left[(1+\nu)\ln{\frac{r_o}{r}}-1\right]\\
    M_t&=\frac{P}{4\pi}\left[(1+\nu)\ln{\frac{r_o}{r}}-\nu\right]		
	\end{split}
\end{equation}

\begin{table}[ht]
\centering
\begin{tabular}{l|ll|ll|}
\cline{2-5}
                              & \multicolumn{2}{l|}{Simple Support}  & \multicolumn{2}{l|}{Fixed Support}   \\ \hline
\multicolumn{1}{|l|}{$\nu$}   & \multicolumn{1}{l|}{Concave} & Convex & \multicolumn{1}{l|}{Concave} & Convex \\ \hline
\multicolumn{1}{|l|}{0.0}     & \multicolumn{1}{l|}{0.9975}  & 0.9982 & \multicolumn{1}{l|}{0.9847}  & 0.9858 \\ \hline
\multicolumn{1}{|l|}{0.3}     & \multicolumn{1}{l|}{0.9928}  & 0.9932 & \multicolumn{1}{l|}{0.9751}  & 0.9749 \\ \hline
\multicolumn{1}{|l|}{0.49}    & \multicolumn{1}{l|}{0.9865}  & 0.9875 & \multicolumn{1}{l|}{0.9638}  & 0.9641 \\ \hline
\multicolumn{1}{|l|}{0.499}   & \multicolumn{1}{l|}{0.9860}  & 0.9870 & \multicolumn{1}{l|}{0.9627}  & 0.9631 \\ \hline
\multicolumn{1}{|l|}{0.49999} & \multicolumn{1}{l|}{0.9863}  & 0.9872 & \multicolumn{1}{l|}{0.9624}  & 0.9624 \\ \hline
\end{tabular}
\caption{Circular Plate Normalized Center Deflection Results (9000 elements), $u_z/u_{theoretical}$}
\label{numplateptloadTab}
\end{table}

\subsection{Elasto-Plastic Cylinder Under Internal Pressure} 
\begin{figure}
\centering
\mbox{
\subfigure[]{\epsfig{file=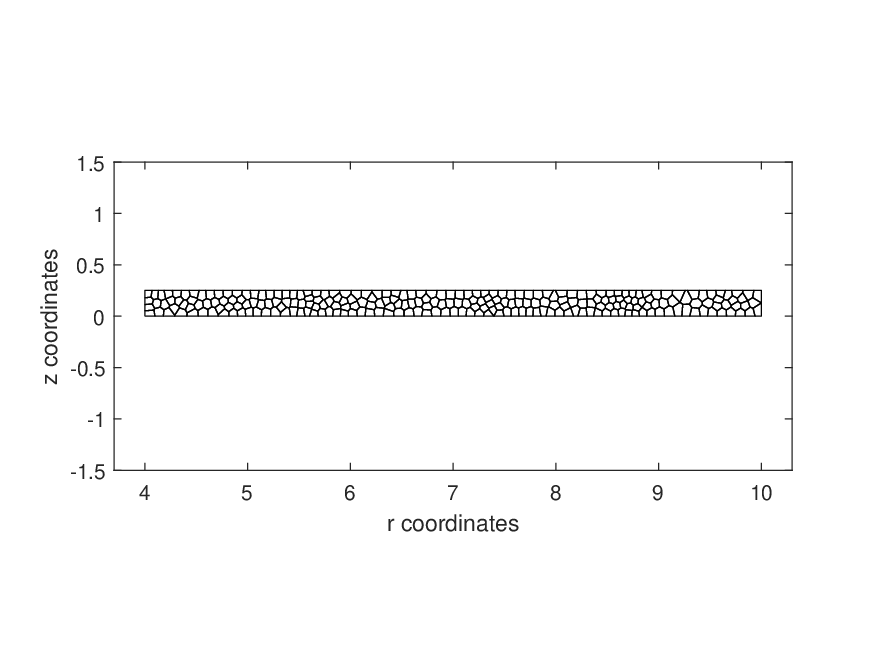,width=0.48\textwidth}}
\subfigure[]{\epsfig{file=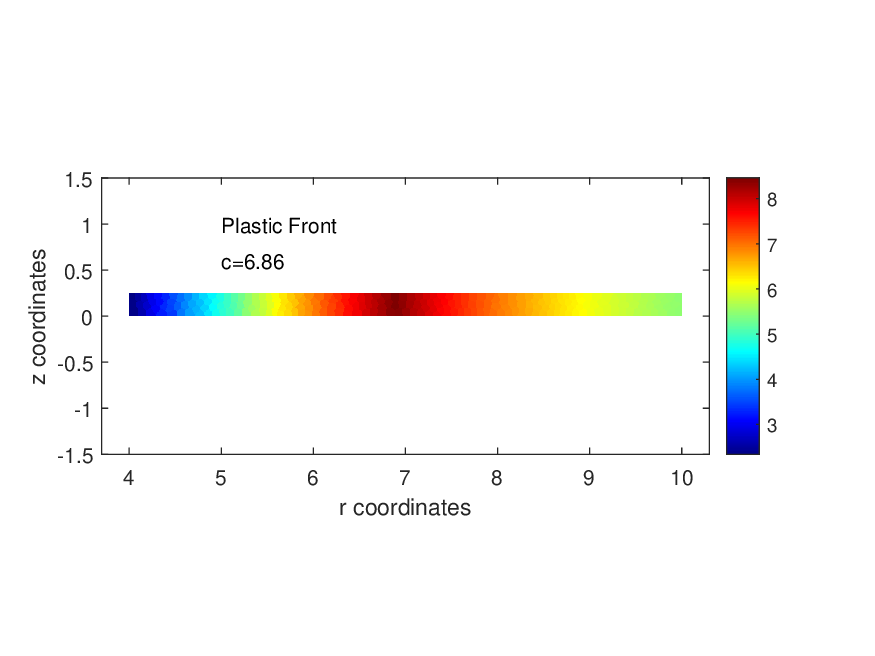,width=0.48\textwidth}}
}

%\vspace*{-0.2in}
\mbox{
\subfigure[]{\epsfig{file=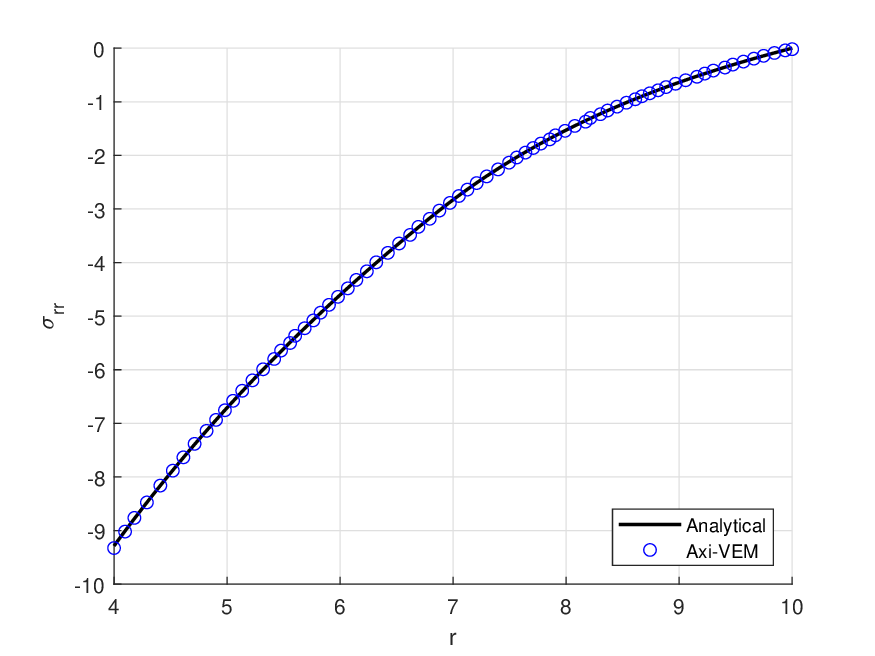,width=0.48\textwidth}}
\subfigure[]{\epsfig{file=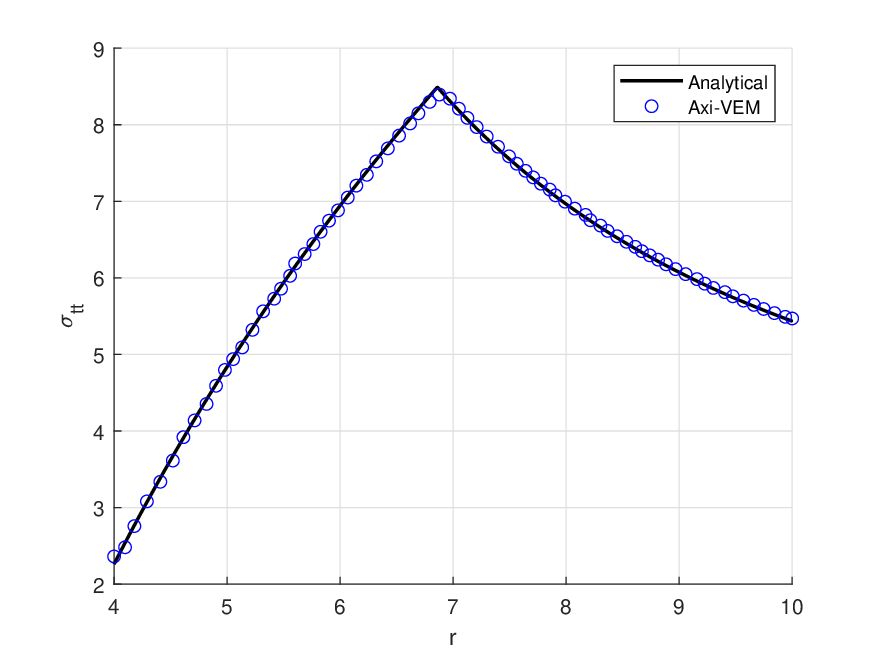,width=0.48\textwidth}}
}

%\vspace*{-0.2in}
%\mbox{
%\subfigure[]{\epsfig{file=./figs/sigmaxring.eps,width=0.48\textwidth}}
%\subfigure[]{\epsfig{file=./figs/sigmayring.eps,width=0.48\textwidth}}
%}

\vspace*{-0.1in}
\caption{Elasto-Plastic Cylinder Under Internal Pressure ($p_{max}=9.29$):  (a)~Axisymmetric cylinder wall mesh (250 convex elements); (b)~Tangential Stress, $\sigma_{tt}$, contour plot; (c)~Radial stress, $\sigma_{rr}$, versus radial coordinate, $r$; and (d)~Tangential stress, $\sigma_{tt}$, versus radial coordinate, $r$.}\label{numericalplasticcyl}
\end{figure}

For the condition of plane strain an elasto-plastic cylinder with 250 convex elements (Figure~\ref{numericalplasticcyl}a) is loaded by an internal pressure, $p$.  The cylinder has inner radius $a=4$, and outer radius, $b=10$.  The modulus of elasticity $E_Y=1000$, Poisson's ratio is $\nu=0.3$, the yield stress is $\sigma_{yield}=10$, and for the case of perfect plasticity the hardening modulus is $E_h=0$.  When a cylinder begins to yield, the `plastic front' location $c$ progresses from the inner radius $a$ toward the outer radius $b$.  Hence, for certain pressures, $a<c<b$.  For the case of pressure, $p=9.29$, the plastic front is located at $c=6.86$, which is in good agreement with the tangential stress results of Figures~\ref{numericalplasticcyl}b and d.   Radial stresses, $\sigma_{rr}$ versus radial coordinate, $r$, are in good agreement with analytical results as shown in Figure~\ref{numericalplasticcyl}c.

The theoretical solution of an elasto-plastic cylinder is given by Hill \cite{Hill}.  Consider the pressure level associated with Figures~\ref{numericalplasticcyl}bcd where the cylinder is yielded and the location of the plastic front, $c$, is such that $a<c<b$.  For the elastic region, where $r \ge c$, the theoretical solution is
\begin{equation}
	\begin{split}
		\sigma_{rr}&=\frac{2\sigma_{yield}}{\sqrt{3}}\frac{c^2}{2b^2}\left(1-\frac{b^2}{r^2}\right)\\
    \sigma_{tt}&=\frac{2\sigma_{yield}}{\sqrt{3}}\frac{c^2}{2b^2}\left(1+\frac{b^2}{r^2}\right).	
	\end{split}
\end{equation}
For the plastic region, where $r<c$, the theoretical solution is
\begin{equation}
	\begin{split}
		\sigma_{rr}&=\frac{2\sigma_{yield}}{\sqrt{3}}\left[-0.5-\ln{\frac{c}{r}}+\frac{c^2}{2b^2}\right]\\
    \sigma_{tt}&=\frac{2\sigma_{yield}}{\sqrt{3}}\left[0.5-\ln{\frac{c}{r}}+\frac{c^2}{2b^2}\right].
	\end{split}
\end{equation}
For a given pressure (a function of c), the plastic front location, $c$, is found numerically from the relation
\begin{equation}
	p(c)=\frac{2\sigma_{yield}}{\sqrt{3}}\left[\ln{\frac{c}{a}}+\frac{1}{2}\left(1-\frac{c^2}{b^2}\right)\right].
\end{equation}

\subsection{Elasto-Plastic Sphere Under Internal Pressure}\label{subsec:plastsphere}
\begin{figure}
\centering
\mbox{
\subfigure[]{\epsfig{file=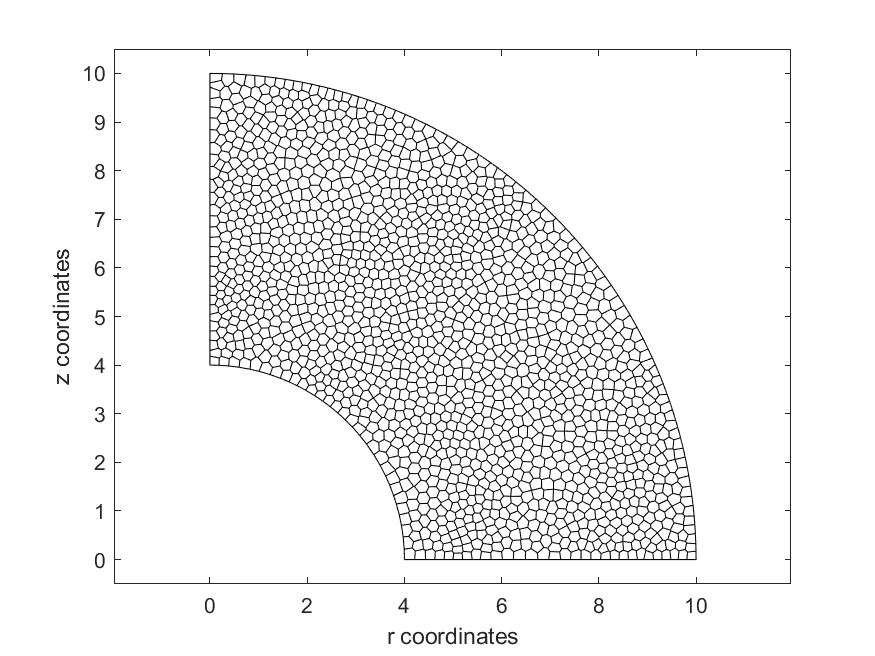,width=0.48\textwidth}}
\subfigure[]{\epsfig{file=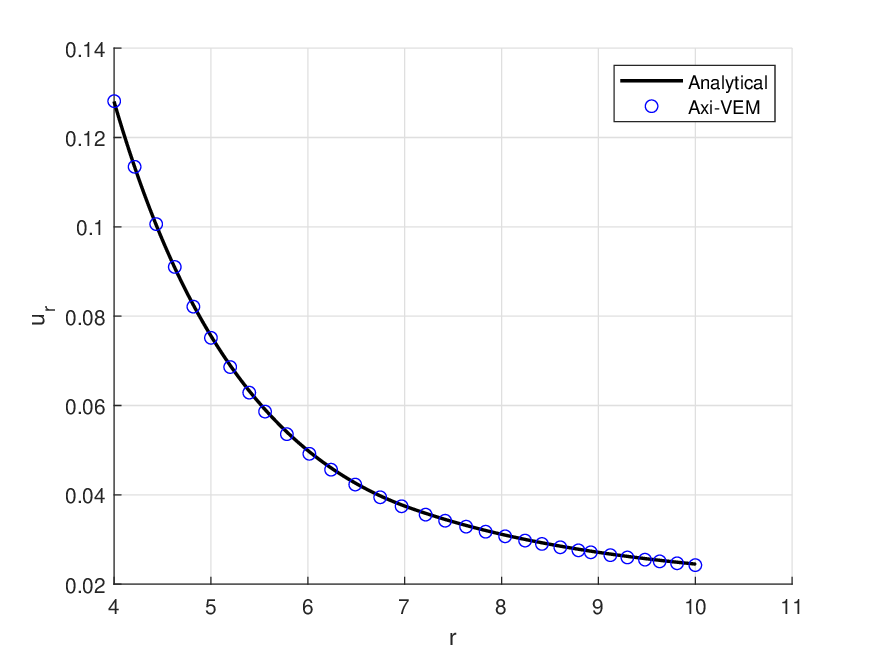,width=0.48\textwidth}}
}

%\vspace*{-0.2in}
\mbox{
\subfigure[]{\epsfig{file=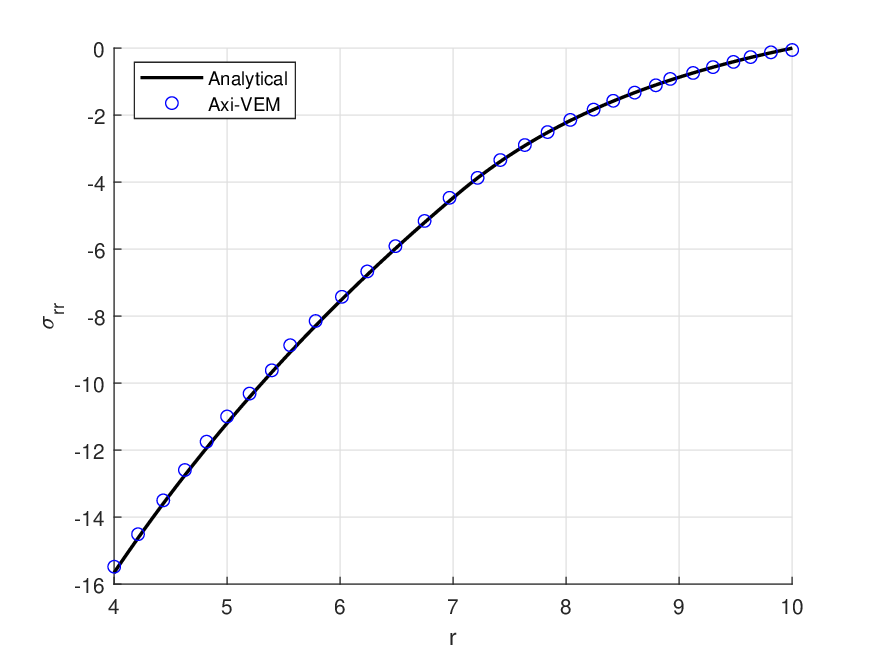,width=0.48\textwidth}}
\subfigure[]{\epsfig{file=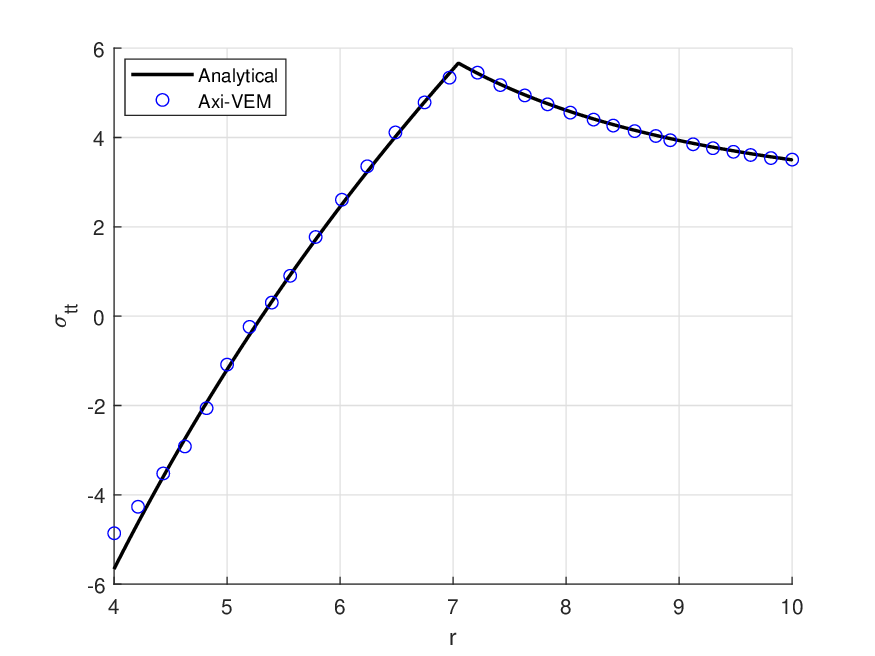,width=0.48\textwidth}}
}

%\vspace*{-0.2in}
\mbox{
\subfigure[]{\epsfig{file=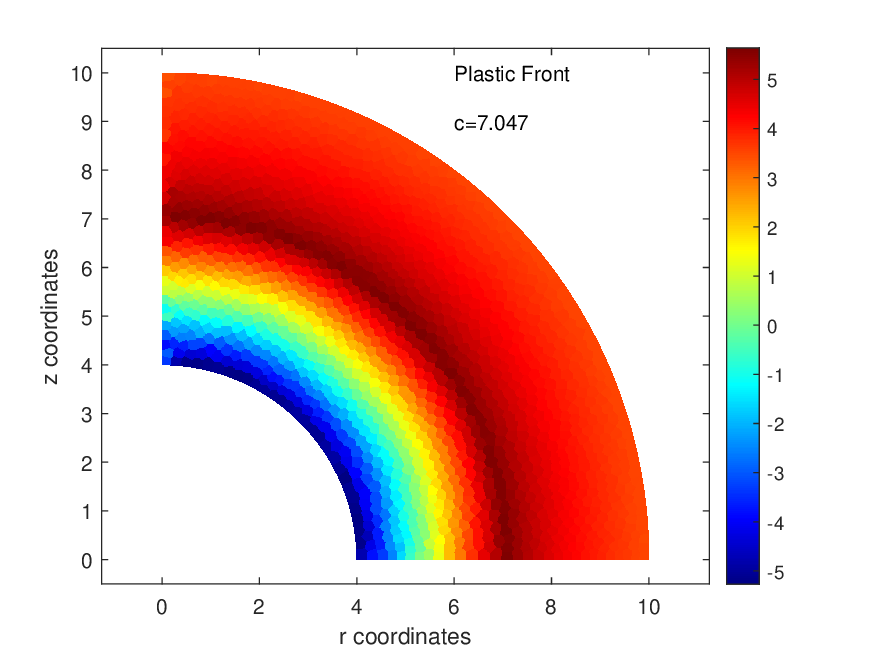,width=0.48\textwidth}}
\subfigure[]{\epsfig{file=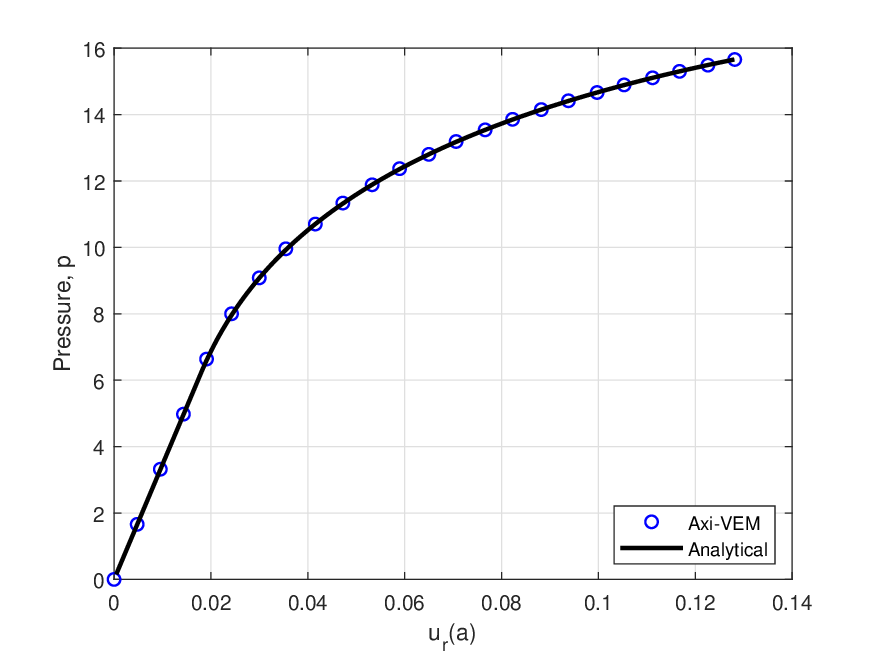,width=0.48\textwidth}}
}

\vspace*{-0.1in}
\caption{Elasto-Plastic Sphere Under Internal Pressure ($p_{max}=15.66$):  (a)~Top of axisymmetric sphere wall mesh (1600 convex elements); (b)~Radial displacement, $u_r$, versus radial coordinate, $r$; (c)~Radial stress, $\sigma_{rr}$, versus radial coordinate, $r$; (d)~Tangential stress, $\sigma_{tt}$, versus radial coordinate, $r$; (e)~Tangential stress, $\sigma_{tt}$, contour plot with plastic front visible at radial distance $c=7.047$; (f) pressure, $p$, versus inner radius wall displacement, $u_r(a)$.}\label{numericalplasticsphere}
\end{figure}

The top half of an elasto-plastic sphere is discretized with 1600 convex elements in Figure~\ref{numericalplasticsphere}a.  The discretization is provided with vertical supports for all nodes at $z=0$ and horizontal supports at all nodes $r=0$.  For the elastic-perfectly plastic sphere material, the modulus of elasticity $E_Y=1000$, Poisson's ratio $\nu=0.3$, yield stress $\sigma_{yield}=10$, and linear hardening modulus $E_h=0$.  The inner and outer radii of the sphere are $a=4$ and $b=10$, respectively.  A nonlinear analysis for the sphere up to a maximum pressure of $p=15.66$ is conducted.  At the maximum pressure Figures~\ref{numericalplasticsphere}bcd show plots for $u_r$, $\sigma_{rr}$, and $\sigma_{tt}$ versus radial coordinate, $r$.  A contour plot of tangential stresses is shown in Figure~\ref{numericalplasticsphere}e.  For the given maximum pressure, $p$, the location of the plastic front, $c=7.047$, is clearly visible in the contour plot and is in agreement with the predicted analytical value (also visible in plot (d)). For the nonlinear analysis a plot of pressure versus internal surface radial displacement $u_r(a)$ is shown in Figure~\ref{numericalplasticsphere}f.  The plot is in excellent agreement with the analytically predicted pressure displacement curve.

The theoretical solution of an elasto-plastic sphere is given by Hill \cite{Hill}.  Consider the pressure level associated with Figures~\ref{numericalplasticsphere}bcde where the sphere is yielded and the location of the plastic front, $c$, is such that $a<c<b$.  For the elastic region, where $r \ge c$, the theoretical solution is
\begin{equation}
	\begin{split}
    u_r&=\frac{2c^3\sigma_{yield}}{3E_Yb^3}\left[(1-2\nu)r+(1+\nu)\frac{b^3}{2r^2}\right]\\
    \sigma_{rr}&=-\frac{2c^3\sigma_{yield}}{3b^3}\left(\frac{b^3}{r^3}-1\right)\\        
    \sigma_{tt}&=\frac{2c^3\sigma_{yield}}{3b^3}\left(\frac{b^3}{2r^3}+1\right).
	\end{split}
\end{equation}
For the plastic region, where $r<c$, the theoretical solution is
\begin{equation}
	\begin{split}
    u_r&=\frac{r\sigma_{yield}}{E_Y}\left[(1-\nu)\frac{c^3}{r^3}-\frac{2}{3}(1-2\nu)\left(1+3\ln{\frac{c}{r}}-\frac{c^3}{b^3}\right)\right]\\   
    \sigma_{rr}&=-\frac{2\sigma_{yield}}{3}\left[1+3\ln{\frac{c}{r}}-\frac{c^3}{b^3}\right]\\    
    \sigma_{tt}&=\frac{2\sigma_{yield}}{3}\left[\frac{1}{2}-3\ln{\frac{c}{r}}+\frac{c^3}{b^3}\right].    
	\end{split}
\end{equation}
For a given pressure (a function of c), the plastic front location, $c$, is found numerically from the relation
\begin{equation}
	p(c)=\frac{2\sigma_{yield}}{3}\left[1+3\ln{\frac{c}{a}}-\frac{c^3}{b^3}\right].
\end{equation}

\subsection{Elasto-Plastic Circular Plate Under Uniform Pressure}\label{subsec:plastcircularplate}
\begin{figure}
\centering
\mbox{
\subfigure[]{\epsfig{file=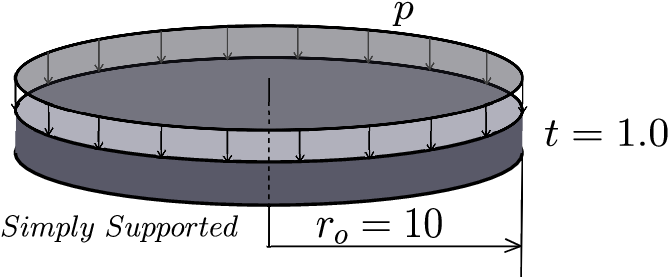,width=0.48\textwidth,clip=}}
\subfigure[]{\epsfig{file=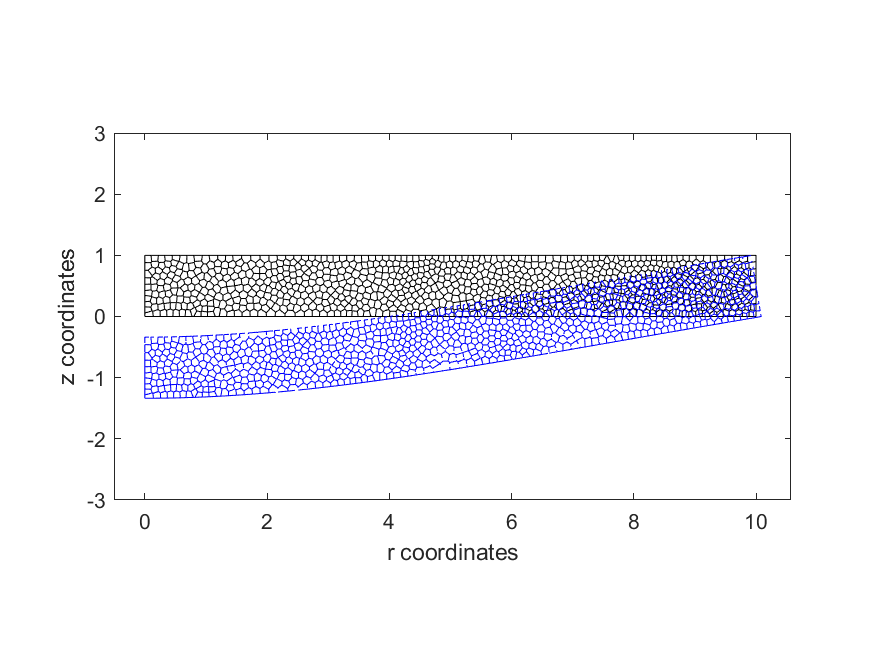,width=0.48\textwidth}}
}

%\vspace*{-0.2in}
\mbox{
\subfigure[]{\epsfig{file=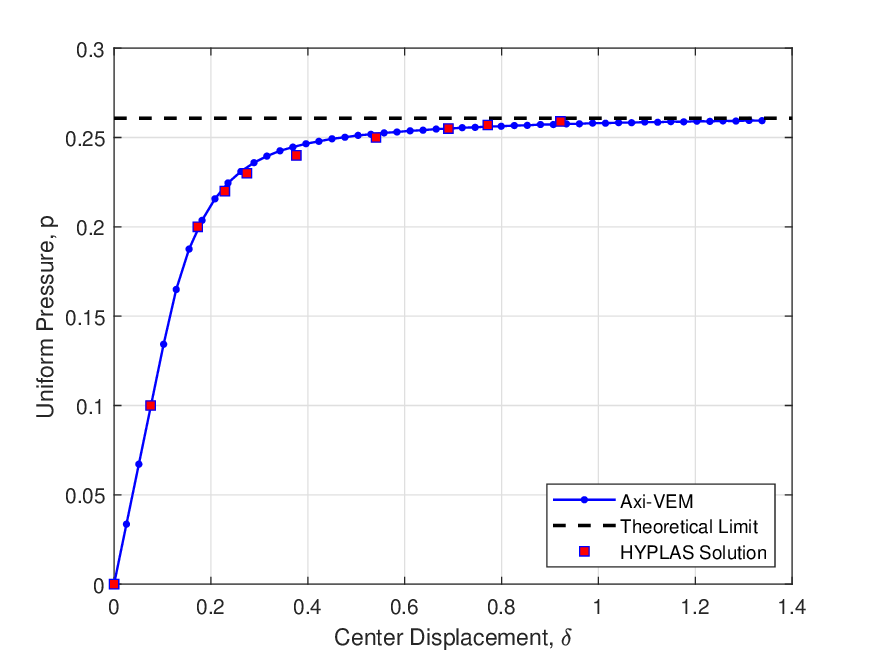,width=0.48\textwidth}}
\subfigure[]{\epsfig{file=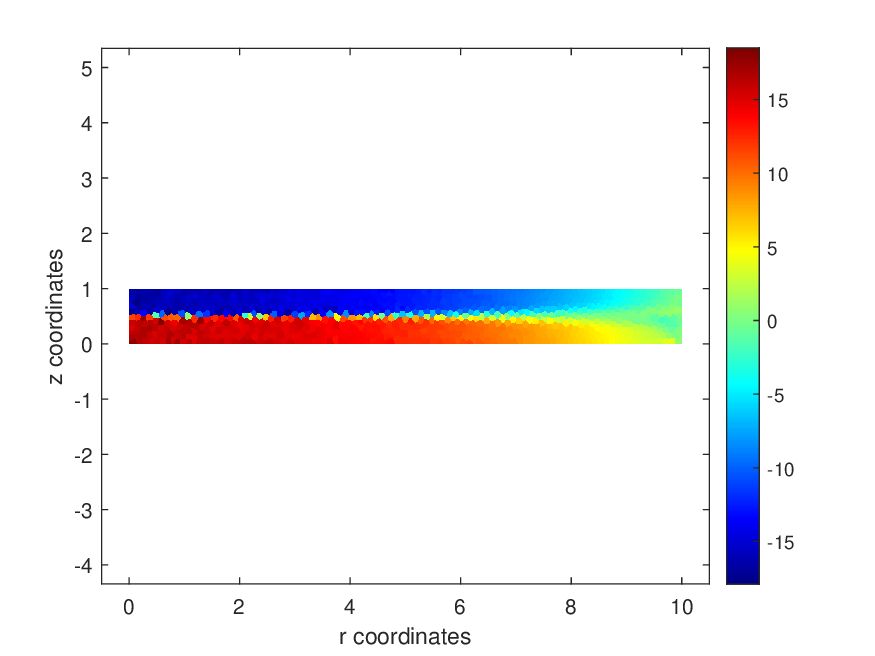,width=0.48\textwidth}}
}

%\vspace*{-0.2in}
%\mbox{
%\subfigure[]{\epsfig{file=./figs/sigmaxring.eps,width=0.48\textwidth}}
%\subfigure[]{\epsfig{file=./figs/sigmayring.eps,width=0.48\textwidth}}
%}

\vspace*{-0.1in}
\caption{Elasto-Plastic Plate Under Uniform Pressure:  (a)~Plate model uniform loading problem; (b)~Undeformed and deflected shape (900 convex elements); (c)~Pressure versus center displacement; and (d)~Contour plot of radial stresses, $\sigma_{rr}$, at maximum pressure of $p=0.2594$.}\label{numericalplasticplate}
\end{figure}

An elastic perfectly plastic simply supported circular plate is modeled with 900 convex polygonal elements.  For the plate material, the modulus of elasticity $E_Y=10000$, Poisson's ratio $\nu=0.24$, yield stress $\sigma_{yield}=16$, and linear hardening modulus $E_h=0$.  The plate of radius $r_o=10$ and thickness $t=1$ is loaded downward with a uniform pressure.  The relevant geometry, loading, discretization, and deflected shape are shown in Figures~\ref{numericalplasticplate}a and b.  A nonlinear plot of pressure versus center displacement is provided in Figure~\ref{numericalplasticplate}c.  The analysis limit pressure, $p=0.2594$, is in very good agreement with the theoretical limit pressure \cite{Skrzypek} $p=0.2609$.  The nonlinear plot is substantially similar to finite element results achieved by the computer program HYPLAS developed by de~Souza~Neto~et~al \cite{deSouzaNeto}. The progression of yielding along the radial direction of the plate is evident in the contour plot of radial stresses, $\sigma_{rr}$, versus radial distance, $r$, shown in Figure~\ref{numericalplasticplate}d.  

% Section 9
\section{Conclusions}\label{sec:conclusions}
In this paper, a first order axisymmetric VEM formulation for elasticity and elasto-plasticity is presented.  This is accomplished by extending a 2D VEM formulation by adding an appropriately constructed row of mean value coordinates shape functions to the strain displacement matrix.  By doing this, tangential strains, needed for axisymmetric problems, are included.  The formulation is limited to small strains for elastic and plastic problems.  Capabilities of the formulation are illustrated by solving example problems and comparing them to theoretical results or the results obtained by finite element solutions.  For both elasticity and plasticity problems the formulation is able to easily obtain results that match the theoretical results.  Several examples also show that near incompressibility is solved without volumetric locking.  For the last problem presented, the axisymmetric VEM formulation with plasticity is able to produce results that are in very good agreement with a finite element solution.  It is evident from the example problems that axisymmetric VEM is a viable method for solving axisymmetric solid mechanics problems.

All numerical methods have limitations.  The formulation herein is no different and has typical limitations observed as follows:  (a) need for suitable incremental step size in nonlinear analysis, particularly for plasticity problems, (b) need for adequate mesh refinement, (c) need for a robust polygon meshing scheme. 

Suggested future work includes the following:  (a) investigate the effect of stabilization free VEM, (b) incorporate finite strains.

\section*{Acknowledgements}
\ack{
LLY acknowledges the research support of Walla Walla University.
Helpful discussions with Alejandro Ortiz-Bernardin are also gratefully acknowledged.
}

\section*{CONFLICT OF INTEREST STATEMENT}The authors declare no potential conflict of interests.

\section*{Data Availability Statement}The data that support the findings of this study are available from the corresponding author upon reasonable request.

\bibliography{bibAxisymmVEM}

\begin{thebibliography}{10}
\providecommand \doibase [0]{http://dx.doi.org/}%

\bibitem{Beirao}
\textrm{Beir\~ao da Veiga} L, Brezzi F, Cangiani A, Manzini G, Marini LD, Russo
  A. Basic Principles of Virtual Element Methods. {\it Mathematical Models and
  Methods in Applied Sciences} 2013\string; 23(1)\string: 199--214.

\bibitem{Beirao2}
\textrm{Beir\~ao da Veiga} L, Brezzi F, Marini LD, Russo A. The Hitchhiker's
  Guide to the Virtual Element Method. {\it Mathematical Models and Methods in
  Applied Sciences} 2014\string; 24(8)\string: 1541--1573.

\bibitem{Beirao:2013:VEL}
\textrm{Beir\~ao da Veiga} L, Brezzi F, Marini D. Virtual elements for linear
  elasticity problems. {\it SIAM J Numer Analysis} 2013\string; 51(12)\string:
  794--812.

\bibitem{Artioli}
Artioli E, \textrm{Beir\~ao da Veiga} L, Lovadina C, Sacco E. {Arbitrary order
  2D virtual elements for polygonal meshes: Part I, elastic problem}. {\it
  Computational Mechanics} 2017\string; 60(3)\string: 355--377.

\bibitem{TaylorVEMplast}
Taylor RL, Artioli E. VEM for Inelastic Solids. In:  O\~{n}ate E. \kern-2pt,
  ed. {\it Advances in Computational Plasticity, Computational Methods in
  Applied Sciences 46\textrm{. }}Switzerland: Springer.  2018 (pp. 381--394).

\bibitem{Artioli2}
Artioli E, \textrm{Beir\~ao da Veiga} L, Lovadina C, Sacco E. {Arbitrary order
  2D virtual elements for polygonal meshes: Part II, inelastic problem}. {\it
  Computational Mechanics} 2017\string; 60(6)\string: 643--657.

\bibitem{Beirao3}
\textrm{Beir\~ao da Veiga} L, Lovadina C, Mora D. A Virtual Element Method for
  elastic and inelastic problems on polytope meshes. {\it Computer Methods in
  Applied Mechanics and Engineering} 2015\string; 295(10)\string: 327--346.

\bibitem{mengolini}
Mengolini M, Benedetto MF, Arag\'{o}n AM. An engineering perspective to the
  virtual element method and its interplay with the standard finite element
  method. {\it Computer Methods in Applied Mechanics and Engineering}
  2019\string; 350\string: 995--1023.

\bibitem{yawVEM}
Yaw LL. Introduction to the Virtual Element Method for 2D Elasticity. {\it
  arXiv: 2301.11928v2 [math.NA]} 2023.

\bibitem{suku:elastodyn}
Sukumar N, Tupek MR. Virtual elements on agglomerated finite elements to
  increase the critical time step in elastodynamic simulations. {\it
  International Journal for Numerical Methods in Engineering} 2022\string;
  123\string: 4702--4725.

\bibitem{Flanagan:1981:AUS}
Flanagan D, Belytschko T. A uniform strain hexahedron and quadrilateral with
  orthogonal hourglass control. {\it International Journal for Numerical
  Methods in Engineering} 1981\string; 17\string: 679--706.

\bibitem{Russo}
Russo A, Sukumar N. Quantitative study of the stabilization parameter in the
  virtual element method. {\it arXiv:2304.00063 [math.NA]} 2023.

\bibitem{Cangiani:2015:HSV}
Cangiani A, Manzini G, Russo A, Sukumar N. Hourglass stabilization and the
  virtual element method. {\it International Journal for Numerical Methods in
  Engineering} 2015\string; 102\string: 404-436.

\bibitem{Hughes}
Hughes TJR. {\it The Finite Element Method - Linear Static and Dynamic Finite
  Element Analysis}.
\newblock Mineola, NY: Dover.
\newblock 1st~ed. 2000.

\bibitem{park}
Park K, Chi H, Paulino G. {B-bar virtual element method for nearly
  incompressible and compressible materials}. {\it Meccanica} 2021\string;
  56\string: 1423--1439.

\bibitem{tb:meshless}
Belytschko T, Krongauz Y, Organ D, Fleming M, Krysl P. Meshless Methods: An
  Overview and Recent Developments. {\it Computer Methods in Applied Mechanics
  and Engineering} 1996\string; 139\string: 3--47.

\bibitem{suku:maxent}
Sukumar N. Construction of polygonal interpolants: A maximum entropy approach.
  {\it International Journal for Numerical Methods in Engineering} 2004\string;
  61(12)\string: 2159--2181.

\bibitem{suku:maxentreview}
Sukumar N, Wright RW. Overview and construction of meshfree basis functions:
  {From} moving least squares to entropy approximants. {\it International
  Journal for Numerical Methods in Engineering} 2007\string; 70(2)\string:
  181--205.

\bibitem{Floater}
Floater MS. Generalized barycentric coordinates and applications. {\it Acta
  Numerica} 2015\string; 24\string: 161--214.

\bibitem{Simo3}
Simo JC, Taylor RL. Consistent Tangent Operators For Rate-Independent
  Elastoplasticity. {\it Computer Methods in Applied Mechanics and Engineering}
  1985\string; 48\string: 101--118.

\bibitem{Simo2}
Simo JC, Hughes TJR. {\it Computational Inelasticity}.
\newblock New York: Springer-Verlag .
\newblock 1998.

\bibitem{Crisfieldv1}
Crisfield MA. {\it Non-linear Finite Element Analysis of Solids and Structures
  -- Vol 1}.
\newblock Chichester, England: John Wiley \& Sons Ltd. .
\newblock 1991.

\bibitem{talischi}
Talischi C, Paulino GH, Pereira A, Menezes IFM. PolyMesher: a general-purpose
  mesh generator for polygonal elements written in Matlab. {\it Structural and
  Multidisciplinary Optimization} 2012\string; 45\string: 309--328.

\bibitem{Chen3}
Chen A, Sukumar N. Stress-hybrid virtual element method on quadrilateral meshes
  for compressible and nearly-incompressible linear elasticity. {\it
  International Journal for Numerical Methods in Engineering} November 1, 2023.

\bibitem{Timoshenko}
Timoshenko SP, Goodier JN. {\it Theory of Elasticity}.
\newblock New York: McGraw-Hill.
\newblock 2nd~ed. 1951.

\bibitem{Timoshenko2}
Timoshenko SP, Woinowsky-Krieger S. {\it Theory of Plates and Shells}.
\newblock New York: McGraw-Hill.
\newblock 2nd~ed. 1959.

\bibitem{Ugural}
Ugural AC. {\it Stresses in Plates and Shells}.
\newblock New York: McGraw-Hill .
\newblock 1981.

\bibitem{Hill}
Hill R. {\it The mathematical theory of plasticity}.
\newblock Oxford, UK: Clarendon Press.
\newblock 1st~ed. 1950.

\bibitem{Skrzypek}
Skrzypek JJ, Hetnarski RB. {\it Plasticity and Creep}.
\newblock Boca Raton, FL: CRC Press .
\newblock 1993.

\bibitem{deSouzaNeto}
{de Souza Neto} EA, Peri\'c D, Owen DRJ. {\it Computational Methods for
  Plasticity}.
\newblock Chichester, UK: John Wiley \& Sons Ltd. .
\newblock 2008.

\end{thebibliography}

\end{document}